%% file: isaac9.tex
\input  amstex
\input amsppt.sty
\magnification1200
\vsize=23.5truecm
\hsize=16.5truecm
\NoBlackBoxes

\input gmacro2.tex

\document
\topmatter
\title
{Fractional-order operators: Boundary problems, heat equations
} 
\endtitle
\author Gerd Grubb \endauthor
\affil
{Department of Mathematical Sciences, Copenhagen University,
Universitetsparken 5, DK-2100 Copenhagen, Denmark.
E-mail {\tt grubb\@math.ku.dk} \footnote{Expanded  
version
of a lecture given at the ISAAC Conference August 2017,
V\"axj\"o{} University, Sweden. \qquad}}\endaffil
\rightheadtext{Fractional-order operators}
\abstract
The first half of this work  gives a survey of the fractional Laplacian (and
related operators), its restricted Dirichlet realization on a bounded
domain, and its nonhomogeneous local boundary conditions, as treated by
pseudodifferential methods. The second half takes
up the associated heat equation with homogeneous Dirichlet
condition. Here we recall recently shown sharp results on interior
regularity and  on $L_p$-estimates up to the
boundary, as well as recent 
H\"older estimates. This is supplied with new higher regularity estimates in
$L_2$-spaces using a technique of Lions and Magenes, and
higher $L_p$-regularity estimates (with arbitrarily high H\"older
estimates in the time-parameter) based on a general result of Amann.
Moreover, it is shown that an improvement to spatial $C^\infty $-regularity
at the boundary is not  in general possible.
\endabstract

\endtopmatter

\head 0. Introduction \endhead

This work is partly a survey of known results for the  fractional
Laplacian and its generalizations, with emphasis on pseudodifferential
methods and local boundary
conditions. Partly it brings new results for the associated heat
equation.

There is an extensive theory for boundary value
problems and evolution problems for elliptic {\bf differential} operators,
developed through many years and including nonlinear problems, and
problems with data of low smoothness. Boundary and evolution problems for
{\bf pseudodifferential} operators, such as fractional powers of the
Laplacian,
have been studied far less, and pose severe difficulties since the
operators
are {\bf nonlocal}.

The presentation here deals with linear questions, since this is
the basic knowledge one needs in any case. Our main purpose is to
explain the application of pseudodifferential methods (with the fractional Laplacian as a
prominent example). The
boundary value theory has been established only in  recent years.
Plan of the paper:
\roster
\item Fractional-order operators.
\item  Homogeneous Dirichlet problems on a subset of ${\Bbb R}^n$.
\item Nonhomogeneous boundary value problems.
\item Heat equations.
\endroster

\example{Remark}
There are other strategies currently in use, such as
methods for singular integral operators in
probability theory and potential theory
(cf.\ e.g.\
\cite{BBC03,CD14,CS98,CT04,
FK13,FR17,J02,JX15,K97,R16,RS14,RS14a,RS15,RSV17,RV18}), methods for embedding the problem in a
degenerate elliptic differential operator situation 
(\cite{CS07} and many subsequent studies, e.g.\ \cite{CSS08}), and  exact
calculations in polar coordinates for the ball (cf.\ e.g.\ \cite{AJS18,DKK17,DG17,ZG16}). 
Each of the methods allow
different types of generalizations of the fractional Laplacian, and
solve a variety of problems, primarily in low-regularity spaces. It is
perhaps surprising that the methods from the calculus of pseudodifferential
operators have only entered the modern studies in the field in the last few years.
\endexample

The new regularity results in Section 4 on heat
 problems, reaching beyond the
recent works
\cite{FR17,G18,RV18}, are: Theorem 4.2 and its corollary giving a limitation
on high spatial regularity,
Theorems 4.6--8 on estimates in 
$L_2$-related spaces for $x$-dependent operators, and
Theorems 4.14--19 on high estimates in H\"older spaces with respect to time,
valued in $L_2$- or $L_p$-related spaces in $x$ (including a
 H\"older-related space in $x$
as a limiting case).

\head 1. Fractional-order operators \endhead

The fractional Laplacian $P=(-\Delta )^a$ on ${\Bbb R}^n$, $0<a<1$, 
has linear and nonlinear 
applications in  
mathematical physics and
differential geometry, and in probability and finance. (See e.\ g.\ Frank-Geisinger \cite{FG16}, Boulenger-Himmelsbach-Lenzmann \cite{BHL16},  Gonzales-Mazzeo-Sire
\cite{GMS12},
Monard-Nickl-Paternain \cite{MNP18}, Kulczycki \cite{K97}, Chen-Song \cite{CS98},  Jakubowski \cite{J02}, Bogdan-Burdzy-Chen
\cite{BBC03}, Applebaum \cite{A04}, Cont-Tankov \cite{CT04}, and
their references.)

The interest in  probability and finance comes from the fact that 
$-P$ generates a semigroup $e^{-tP}$ which is a stable
L\'evy process. Here $P$ is viewed as a {\bf  singular integral operator:} 
$$
(-\Delta )^au(x)=c_{n,a}PV\int_{{\Bbb
    R}^n}\frac{u(x)-u(x+y)}{|y|^{n+2a}}\,dy;
\tag 1.1
$$
 more general stable L\'evy processes arise from operators
$$
Pu(x)=PV\int_{{\Bbb
    R}^n}(u(x)-u(x+y))K(y)\,dy,\quad K(y)=\frac{K(y/|y|)}{|y|^{n+2a}},\tag 1.2
$$
where the homogeneous kernel function $K(y)$ is locally integrable, 
 {\bf  positive}, and {\bf
  even}: $K(-y)=K(y)$, on ${\Bbb
  R}^n\setminus \{0\}$. (1.2) can also be generalized to nonhomogeneous
  kernel functions satisfying suitable estimates in terms of $|y|^{-n-2a}$.
Usually, only real functions are considered in probability studies.

$(-\Delta )^a$ can instead be viewed as a {\bf  pseudodifferential operator} ($\psi $do)
of order $2a$:
$$
(-\Delta )^au=\operatorname{Op}(|\xi |^{2a})u=
\Cal F^{-1}(|\xi |^{2a}\Cal F u(\xi )),\tag 1.3
$$
 using the Fourier transform $\Cal F$, defined by $\hat
u(\xi )=\Cal F
u(\xi )= \int_{{\Bbb R}^n}e^{-ix\cdot \xi }u(x)\, dx$ (extended from the
Schwartz space $\Cal S({\Bbb R}^n)$ of rapidly decreasing $C^\infty
$-functions, to the temperate distributions  $\Cal S'({\Bbb R}^n)$). 
  $\Psi $do's are in general defined by
$$
Pu=\operatorname{Op}(p(x,\xi ))u=
\Cal F^{-1}_{\xi\to x}(p(x,\xi )\Cal F u(\xi ));\tag 1.4
$$
note that this theory operates in a context of complex functions (and
distributions). Moreover, (1.4) allows $x$-dependence in the symbol $p(x,\xi )$.

In (1.2), if 
$K\in C^\infty (\Bbb R^n\setminus \{0\})$, the operator it defines is
the same as the operator defined by $p(\xi )=\Cal F
K(y)$ in (1.4); here $p(\xi )$ is homogeneous of degree $2a$,
positive and even.

As a generalization of (1.3), we  consider $x$-dependent classical
$\psi $do's of order $2a\in\rp$, with certain properties. That $P$ is
classical of order $2a$ means that there is an asymptotic expansion of
the symbol $p(x,\xi )$ in a series of terms $p_j(x,\xi )$, $j\in {\Bbb N}_0$, that are
homogeneous in $\xi $ of order $2a-j$ for $|\xi |\ge 1$; the expansion
holds in the sense that $$
|\partial_\xi ^\alpha \partial_x^\beta [p(x,\xi )-
 \sum_{j<M}p_j(x,\xi )]|\le C_{a,\beta ,M}\ang\xi ^{2a-|\alpha |-M}, \text{
 all }\alpha ,\beta \in{\Bbb N}_0^n, M\in{\Bbb
 N}_0;\tag1.5$$
here $\langle\xi \rangle$ stands for $(|\xi |^2+1)^{\frac12}$. To the operators
defined from these symbols by (1.4) one adds the smoothing operators
mapping $\E'({\Bbb R}^n)$ to $C^\infty ({\Bbb R}^n)$ (also called negligible operators).
We assume moreover that  $p$ is {\bf even}, meaning that
$$
p_j(x,-\xi )=(-1)^jp_j(x,\xi ),\text{ all }j, \tag 1.6
$$
and {\bf strongly elliptic}, meaning that for a positive constant $c$,
$$
\operatorname{Re}p_0(x,\xi )\ge c|\xi |^{2a} \text{ for }|\xi |\ge 1.\tag1.7
$$
Then $P=\Op (p(x,\xi ))$ can be shown to have some of the same features as
$(-\Delta )^a$. To sum up, we are assuming, with $a\in{\Bbb R}_+$:

\definition{Hypothesis 1.1} $P$ is a classical pseudodifferential operator of order
$2a$, even and strongly elliptic (cf.\ (1.4), (1.6), (1.7)).
\enddefinition

In part of Section 4, we moreover assume $a<1$. For some results we consider the subset of operators defined as in
(1.2)ff.\ with a kernel function $K(y)$ that is
smooth outside 0: 

\definition{Hypothesis 1.2}
$P$ is as in {\rm (1.2)}, with $K(y)$ positive, homogeneous of degree
$-n-2a$, even, and $C^\infty $ on ${\Bbb R}^{n-1}\setminus \{0\}$. 
\enddefinition

This fits into the $\psi $do formulation, when we 
write $p(\xi )=\Cal F K(y)$ as
$(1-\chi (\xi ))p(\xi )
+\chi (\xi )p(\xi )$, with $\chi \in C_0^\infty ({\Bbb
R}^n)$ and $\chi (\xi )=1$ near 0. Here the first term is a symbol as
under  Hypothesis 1.1, and the operator defined from the second term
maps large spaces of distributions (e.g.\ $\E'({\Bbb R}^n)$) into
$C^\infty $-functions (hence is a negligible operator).

To give an example of an $x$-dependent operator, we can mention that $(-\Delta )^a$ will take the $x$-dependent form
if it undergoes a smooth change of coordinates. As a more general example,  $P$ can be
an operator defined as  $P=A(x,D)^a$,
where $A(x,D)$ is a second-order strongly elliptic differential operator.
Here $P$ is constructed via the resolvent (Seeley \cite{S69}). But of
course, the symbol $p(x,\xi )$ can be taken much more general, not
tied to differential operator considerations.

A difficult aspect of such operators is that they are {\bf
  nonlocal}. This is a well-known feature in the pseudodifferential
  theory, where one can profit from pseudo-locality (namely, $Pu$ is $C^\infty
  $ on the set where $u$ is $C^\infty $). In a different approach,  Caffarelli and Silvestre \cite{CS07} showed that  $(-\Delta )^a$ on $\Bbb R^n$ is the Dirichlet-to-Neumann
operator for a degenerate elliptic {\bf differential} boundary value
problem  on $\Bbb R^n\times \Bbb R_+$;  {\bf local} in dimension $n+1$.
This observation  was then used 
to obtain results by transforming problems for $(-\Delta
)^a$ into problems for  {\bf local} operators in one more variable, e.g.\ in
  \cite{CSS08}.
(However, in some cases where one needs to consider $(-\Delta )^a$ 
over a subset $\Omega \subset{\Bbb R}^n$, the transformation might lead to 
equally difficult problems in the new variables.)

\head{2. Homogeneous Dirichlet problems on a subset of ${\Bbb R}^n $}\endhead

How do we get $P$ to act over $\Omega $? There are several ways to
answer this.
Let us first introduce an appropriate scale of $L_p$-based Sobolev
spaces. 

The standard Sobolev-Slobodetski\u \i{} spaces $W^{s,p}({\Bbb R}^n)$, $1<p<\infty $ and
$s\ge 0$, have a different character according to whether $s$ is
integer or not. Namely, for $s$ integer, they consist of
$L_p$-functions with derivatives in $L_p$ up to order $s$, hence
coincide with the Bessel-potential spaces $H^s_p({\Bbb R}^n)$, defined
for $s\in{\Bbb R}$ by 
$$
H_p^s(\R^n)=\{u\in \SD'({\Bbb R}^n)\mid \F^{-1}(\ang{\xi }^s\hat u)\in
L_p(\R^n)\}.\tag2.1
$$
For noninteger $s$, the $W^{s,p}$-spaces coincide with the Besov spaces, defined e.g.\ 
as follows: For $0<s<2$,
$$
f\in
B^s_p({\Bbb R}^{n})\iff \|f\|_{L_p}^p+ \int_{\R^{2n}}\frac{|f(x)+f(y)-2f((x+y)/2)|^p}{|x+y|^{n+ps}}\,dxdy<\infty ;\tag2.2
$$
and $B^{s+t}_p(\R^n)=(1-\Delta )^{-t/2}B^s_p({\Bbb R}^n)$ for all $t\in{\Bbb R}$.

The Bessel-potential spaces $H^s_p$ are important because they are most
directly related to $L_p$; the Besov spaces $B^s_p$ have other
convenient properties, and are
needed for boundary value problems in an $H^s_p$-context,
 because they are the correct range spaces for
trace maps (both from $H^s_p$ and $B^s_p$-spaces); see e.g.\ the overview in
the introduction to \cite{G90}. For $p=2$, the two scales are
identical, and $p$ is usually omitted. For $p\ne 2$ they are related by strict inclusions:$$
H^s_p\subset B^s_p\text{ when }p>2,\quad H^s_p\supset B^s_p\text{ when }p<2.\tag 2.3
$$
When working with operators of noninteger order, the use of the
$W^{s,p}$-notation can lead to confusion since the definition depends
on the integrality of $s$; moreover, this scale does not
always interpolate well.
In the following, we focus on the Bessel-potential scale $H^s_p$, but much of
what we show is directly generalized to the Besov scale $B^s_p$, and to other
scales (Besov-Triebel-Lizorkin spaces). There is a more general Besov scale $B^s_{p,q}$ (cf.\ e.g.\
Triebel \cite{T78}), where $B^s_{p} $ equals the special case $B^s_{p,p}$.

There is an identification of $H^s_p({\Bbb R}^n)$ with the dual space
of $H^{-s}_{p'}({\Bbb R}^n)$, $1/p+1/p'=1$, in a duality consistent with
the $L_2$-duality, and there is a similar result for the Besov scale.

Let $\Omega $ be a open subset of ${\Bbb R}^n$ (we shall use it with
$C^\infty $-boundary, but much of the following holds under limited
smoothness assumptions). One defines the two
associated scales relative to $\Omega $ (the {\it restricted} resp.\ {\it supported} version):
$$
\aligned
 \overline H_p^s(\Omega )&=r^+H^s_p(\Bbb R^n), \\
\dot H^s_p(\overline\Omega )&=\{u\in H^s_p(\Bbb R^n)\mid
\operatorname{supp}u\subset \overline\Omega  \};
\endaligned\tag 2.4
$$
 here $\operatorname{supp}u$ denotes the support of $u$ (the
 complement of the largest open set where $u$ is zero).
Restriction from $\R^n$ to $\Omega $  is denoted $r^+$,
 extension by zero from $\Omega $ to $\R^n$ is denoted $e^+$ (it is
 sometimes tacitly understood).  Restriction
 from $\comega$ to $\partial\Omega $
 is denoted $\gamma _0$.

When $s>1/p-1$, one can identify  $\dot H^s_p(\overline\Omega )$ with a
subspace of $ \overline H_p^s(\Omega )$, closed if $s-1/p\notin {\Bbb
N}_0$ (equal if $1/p-1<s<1/p$), and with a stronger norm if $s-1/p\in{\Bbb N}_0$.

The space $\dot H^s_p$ is in some texts indicated with a ring, zero or twiddle, as e.g.\
$\overset\circ\to H^s_p$, $H^s_{p,0}$ or $\widetilde H^s_p$. In most current texts, $\ol
H_p^s(\Omega )$ is denoted $H_p^s(\Omega )$ without the overline (that
was introduced along with the notation $\dot H_p$ in \cite{H65,H85}),
but we prefer to use it, since it is
 makes the role of the space more
clear in formulas where
both types occur.

Now let us present some operators associated with $(-\Delta )^a$ on $\Omega $:

\smallskip
{\bf (a)} The {\bf restricted} Dirichlet fractional Laplacian
$P_{\operatorname{Dir}}$. It acts like $P=(-\Delta )^a$,
defined on
functions $u$ that are 0 on ${\Bbb R}^n\setminus \Omega $, and followed by
restriction $r^+$ to $\Omega $:
$$
P_{\operatorname{Dir}}u\text{ equals }r^+Pu\text{ when }\operatorname{supp}u\subset \overline\Omega .\tag2.5
$$
 
In $L_2(\Omega )$, it  is the operator defined variationally from the sesquilinear form
$$
Q_0(u,v)=\tfrac12 c_{n,a} \int_{{\Bbb
R}^{2n}}\frac{(u(x)-u(y))(\bar v(x)-\bar v(y))}{|x-y|^{n+2a}}
\,dxdy,\; D(Q_0)=\dot H^a(\overline\Omega ).\tag2.6
$$
 
{\bf (b)} The {\bf spectral} Dirichlet fractional Laplacian $(-\Delta
_{\operatorname{Dir}})^a$, defined  e.g.\ via eigenfunction expansions
of $-\Delta _{\operatorname{Dir}}$. It does not act like $r^+Pe^+$. It
is not often used in probability applications. (Its regularity
properties in $L_p$-Sobolev spaces are discussed in \cite{G16}, which
gives many references to the literature on it.)

{\bf (c)}  The {\bf regional} fractional Laplacian, defined from the
sesquilinear form
$$
Q_1(u,v)=\tfrac12 c_{n,a} \int_{\Omega \times\Omega }\frac{(u(x)-u(y))(\bar v(x)-\bar v(y))}{|x-y|^{n+2a}}
\,dxdy,\;D(Q_1)=\overline H^a(\Omega ).\tag2.7
$$
 It acts like $r^+Pe^++w$ with a certain correction function $w$.

\smallskip
There are still other operators over $\Omega $ that can be defined from
$P$, e.g.\ representing suitable Neumann problems. A local Neumann condition will be discussed below in
Section 3. We refer to \cite{G16} Sect.\ 6, and its references, for an overview over the various choices.

We shall now focus on {\bf (a)}, where the operator acts like $r^+P$.

The homogeneous Dirichlet problem, for a smooth bounded open set
$\Omega $, is
$$
r^+Pu= f \text{ in }\Omega ,\quad \operatorname{supp}u\subset
\overline\Omega .\tag2.8
$$
As $P$ we take $(-\Delta )^a$, or a more general $\psi $do as in
Hypothesis 1.1 or 1.2.
 
$P_{\operatorname{Dir}}$ in  $L_2(\Omega )$ is the realization of $r^+P$ with domain
$$
D(P_{\operatorname{Dir}})=\{u\in \dot H^a(\overline\Omega )\mid r^+Pu\in L_2(\Omega )\}.\tag2.9
$$
When $P$ satisfies Hypothesis 1.2 (in particular, when  $P=(-\Delta
)^a$), then $P_{\operatorname{Dir}}$ is positive
selfadjoint; for other $P$ it is sectorial, with discrete
spectrum in a sector. What can be said about the regularity of
functions in the
domain?
\smallskip

$\bullet$ Vishik and Eskin showed in the 1960's (see e.g.\ Eskin
\cite{E81]}: 
$$D(P_{\operatorname{Dir}})=\dot
H^{2a}(\overline\Omega )\text{ if }a<\tfrac12,\; D(P_{\operatorname{Dir}})\subset\dot
H^{a+\frac12-\varepsilon }(\overline\Omega )\text{ if
}a\ge\tfrac12.\tag 2.10$$

$\bullet$ Ros-Oton and Serra \cite{RS14} showed in 2014 for $(-\Delta
)^{a}$: 
$$f\in L_\infty
(\Omega ) \implies
u\in d^aC^\alpha (\overline\Omega )\text{ for small }\alpha;\tag2.11$$
 here
$d(x)$ equals $\operatorname{dist}(x,\partial\Omega )$ near
$\partial\Omega $, and $C^\alpha $ is the H\"older space. They improved this later to
$\alpha <a$ ($\alpha =a$ in some cases), and to more general
 $P$ as in (1.2), and lifted the regularity conclusions to $f\in C^\gamma $,
$u/d^a\in C^{a+\gamma }$ for small $\gamma $. For (2.11), $\Omega $
was assumed to be $C^{1,1}$.

$\bullet$ We showed in 2015 \cite{G15}, for $1<p<\infty $ and $\Omega $ smooth:  
$$
\aligned
f\in \overline  H^s_p(\Omega )&\iff
u\in H_p^{a(s+2a)} (\overline\Omega ), \text{ any }s\ge 0,\\
f\in C^\infty (\comega )&\iff
u/d^a\in C^\infty (\overline\Omega )
;\endaligned\tag2.12$$
here $ H_p^{a(s+2a)} (\overline\Omega )$ is a space introduced by
H\"ormander \cite{H65} for $p=2$. E.g.\ when $s=0$,
$$
H_p^{a(2a)} (\overline\Omega )\;\cases=\dot
  H_p^{2a}(\overline\Omega )\text{ if }a<1/p,\\
\subset\dot
  H_p^{2a-\varepsilon }(\overline\Omega )\text{ if }a=1/p,\\
\subset\dot
  H_p^{2a}(\overline\Omega )+d^a \overline H_p^a(\Omega )\text{ if }a>1/p;
\endcases\tag2.13
$$
the spaces will be further explained below. (2.12) has corollaries in H\"older spaces by Sobolev embedding.

\smallskip
The contribution from H\"ormander, accounted for in detail in
\cite{G15}, is in short the following: He
defined the {\bf $\mu $-transmission property} in his book
1985, Sect.\ 18.2: 

\proclaim{ Definition 2.1} A classical $\psi $do  $P$ of order $m$ has the $\mu  $-transmission
  property at 
  $\partial\Omega $,  when
$$
\partial_x^\beta \partial_\xi ^\alpha {p_j}(x,-\nu)=e^{\pi i(m-2\mu  -j-|\alpha | )
}\partial_x^\beta \partial_\xi ^\alpha{p_j}(x,\nu), \tag2.14
$$
for all indices; here $x\in \partial\Omega $, and  $\nu$ denotes the
interior normal at $x$.
\endproclaim

The property was formulated already in a photocopied lecture note fom IAS
Princeton 1965-66 on $\psi $do boundary problems \cite{H65}, handed out to a few
people through the times, including
Boutet de Monvel in 1968, the present author in 1980.  
 
For $P$  of order $2a$ and {\it even} (cf.\ (1.6)), it holds with $\mu
 =a$, for {\it any} smooth subset $\Omega $ ({\it all} normal directions are
 covered when (1.6) holds).

 The case $\mu =0$ is the transmission condition entering in the
 calculus of Boutet de Monvel, described e.g.\ in \cite{B71,G96,S01,G09}.

 Recalling that
$e^+$ denotes extension by zero, let
$$
\Cal E_a(\overline\Omega )=e^+d^aC^\infty (\overline\Omega  ). \tag2.15
$$
Then by \cite{H85}, Th.\ 18.2.18,

 \proclaim{
 Theorem 2.2}   The $a$-transmission property at $\partial\Omega $ is necessary and sufficient in order
that $r^+P$ maps
$\Cal E_a(\overline\Omega )$ into $ C^\infty (\overline\Omega )$.

\endproclaim

Note the importance of $d^a$.
 
The notation in \cite{H85} is slightly different from that in the
notes \cite{H65}, which we have adapted here. 
The  notes  moreover treated solvability questions, with $f\in C^\infty
(\overline\Omega )$ or in $H^s$-spaces. The space
$H^{a(s)}(\overline\Omega )$ was 
introduced. Originally it was defined as ``the functions supported in
$\overline\Omega $ that are mapped into $\overline H^{s-m}(\Omega )$
for any $P$ that is elliptic of order $m$ and has the $a$-transmission
property'', and the
whole effort was to sort this out. We shall now explain the structure
and its implications, for
general $H^s_p$-spaces with $1<p<\infty $.

Introduce first {\bf order-reducing operators of plus/minus type.} For $\Omega ={\Bbb
  R}^n_+$, define for $t\in\Bbb R$:
$$\Xi _\pm^t=\operatorname{Op}((\langle{\xi '}\rangle\pm i\xi_n)^t)\text{ on }\Bbb
R^n.\tag2.16$$
 
The  symbols extend analytically in $\xi _n$ to $\operatorname{Im}\xi
_n\lessgtr 0$. Hence, by the Paley-Wiener theorem, $\Xi _\pm^t$
preserve support in $\overline{\Bbb R}^n_\pm$.
Then for all $s\in\Bbb R$, 
$$
\aligned
\Xi ^t_+&\colon \dot H_p^s (\overline{\Bbb R}^n_+)\overset\sim\to\rightarrow \dot
H_p^{s-t} (\overline{\Bbb R}^n_+),\text{ with inverse } \Xi ^{-t}_+,\\
r^+\Xi ^t_-e^+&\colon \overline H_p^{s} ({\Bbb R}^n_+)\overset\sim\to\rightarrow
\overline H_p^{s-t} ({\Bbb R}^n_+),\text{ with inverse }r^+\Xi
^{-t}_-e^+.
\endaligned\tag2.17
$$
Here the action of $e^+$ on spaces with $s<0$ is understood such that
the operators in the families  $ \Xi ^{t}_+$ and  $ r^+\Xi ^{t}_-e^+$
are adjoints for each $t\in{\Bbb R}$:
$$
\Xi ^t_+\colon \dot H_p^s (\overline{\Bbb R}^n_+)\overset\sim\to\rightarrow \dot
H_p^{s-t} (\overline{\Bbb R}^n_+)\text{ has the adjoint }
r^+\Xi ^t_-e^+\colon \overline H_{p'}^{\,-s+t} ({\Bbb R}^n_+)\overset\sim\to\rightarrow
\overline H_{p'}^{\,-s} ({\Bbb R}^n_+),\tag2.18
$$
with respect to an extension of the duality $\int_{{\Bbb R}^n_+}u\bar v
\, dx $ (more explanation in \cite{G15}, Rem.\ 1.1).

Now define the  {\bf $a$-transmission
  space} over $\Bbb R^n_+$:
$$
H_p^{a(s)} (\overline{\Bbb R}^n_+)=\Xi _+^{-a}e^+\overline
H_p^{s-a} ({\Bbb R}^n_+), \text{ for }s-a>-1/p'.\tag2.19
$$
 
Here $e^+\overline H_p^{s-a} ({\Bbb R}^n_+)$ has a jump
at $x_n=0$ when $s-a>1/p$; this is mapped by $\Xi _+^{-a}$ to a singularity of the type
$x_n^a$. 
 
In fact, we can show:
$$
 H_p^{a(s)}  (\overline{\Bbb R}^n_+)
\cases
=\dot H_p^{s} (\overline{\Bbb R}^n_+)\text{ if }-1/p'<s-a<1/p,\\
\subset \dot
H_p^{s} (\overline{\Bbb R}^n_+)+e^+ x_n^a\overline H_p^{s-a} ({\Bbb R}^n_+)\text{ if }s-a-1/p\in \Bbb
R_+\setminus \Bbb N,
\endcases\tag2.20
$$
with $\dot
H_p^{s} (\overline{\Bbb R}^n_+)$ replaced by $\dot
H_p^{s-\varepsilon } (\overline{\Bbb R}^n_+)$ if $s-a-1/p\in \Bbb N$.

For example, for $1/p<s-a< 1+1/p$, $u\in H_p^{a(s)}(\crnp)$
has the form
$$
u=w+e^+x_n^aK_0\varphi ,\tag2.21
$$
where $w$ and $\varphi $
run through $\dot H_p^{s}(\crnp)$ and $
B_p^{s-a-1/p}({\Bbb R}^{n-1} )$, respectively,
and $K_0 $ is the Poisson operator $K_0\colon \varphi \mapsto \Cal F^{-1}_{\xi
'\to x'}[\hat \varphi (\xi ')r^+e^{-\ang{\xi '}x_n}]$ solving the
standard Dirichlet problem
$$
(-\Delta +1)v=0 \text{ in }\rnp,\quad \gamma _0u=\varphi \text{ on }{\Bbb R}^{n-1}.
$$

The analysis hinges on the following formula for the inverse Fourier transform
of \linebreak
$(\langle{\xi '}\rangle+ i\xi_n)^{-a-1}$, where $e^+r^+x_n^a$ appears:
$$
{\Cal F}^{-1}_{\xi _n\to x_n}(\ang{\xi '} +i\xi _n)^{-a-1}={\Gamma
(a+1)^{-1}} e^+r^+x_n^a e^{-\ang{\xi '} x_n}.
$$

The  generalization  to $\Omega \subset {\Bbb R}^n$ depends on finding
suitable replacements of $\Xi ^t_\pm$. They are a kind of generalized
$\psi $do's (the symbols satify some but not all of the usual symbol
estimates). It was important in \cite{G15} that we could rely on a
truly pseudodifferential version $\Lambda _\pm^{(t)}$ found in \cite{G90}.

The choice $P=(1-\Delta )^a$ on $\Bbb R^n$ with symbol $(1+|\xi
|^2)^a$ serves as a model case with easy explicit calculations. Here one can factorize the
symbol and operator directly:
$$
(1+|\xi |^2)^a=(\langle{\xi '}\rangle- i\xi_n)^a(\langle{\xi
  '}\rangle+ i\xi_n)^a,\quad (1-\Delta )^a=\Xi_- ^a \Xi_+ ^a.\tag2.22
$$
Let us show how to solve the model Dirichlet problem
$$
r^+(1-\Delta )^au=f\text{ on }{\Bbb R}^n_+,\quad \operatorname{supp}u\subset
\overline{\Bbb R}^n_+.\tag2.23
$$
 Say, $f$ is given in 
$\overline
 H_p^t({\Bbb R}^n_+)$ for some $t\ge 0$, and $u$ is a priori assumed
 to lie in $ \dot
 H_p^a(\overline {\Bbb R}^n_+)$.

In view of the factorization (2.22),
$$
r^+(1-\Delta )^au=r^+\Xi ^a_-\Xi ^a_+u=r^+\Xi ^a_-(e^+r^++e^-r^-)\Xi ^a_+u=r^+\Xi ^a_-e^+r^+\Xi ^a_+u,
$$
since $r^-\Xi ^a_+u=0$. ($r^-$ denotes restriction from ${\Bbb R}^n$
to $\rnm$, $e^-$ is extension by zero on  ${\Bbb R}^n\setminus \rnm$.)
By (2.17), the problem
(2.23) is reduced by
composition with $r^+\Xi ^{-a}_-e^+$ to the left to the problem
$$
r^+\Xi ^a_+u=g
,\quad \operatorname{supp}u\subset
\overline{\Bbb R}^n_+,\tag2.24
$$
where $g=r^+\Xi ^{-a}_-e^+f\in \overline H_p^{t+a}({\Bbb R}^n_+)$. Now
there is an important observation, shown in Prop.\ 1.7 in \cite{G15}:

\proclaim{Lemma 2.3} Let $s>a-1/p'$. The mapping $\Xi _+^{-a}e^+$ is a
bijection from $\ol H_p^{s-a}(\rnp)$ to $H_p^{a(s)}(\crnp)$ with
inverse $r^+\Xi _+^a$.
\endproclaim

Then clearly, (2.24) is simply solved uniquely by 
$$
u=\Xi _+^{-a}e^+g.\tag2.25
$$

Inserting the definition of $g$, we can conclude:

\proclaim{Proposition 2.4} The problem {\rm (2.23)} with 
$f$ is given in 
$\overline
 H_p^t({\Bbb R}^n_+)$ for some $t\ge 0$, and $u$ sought in $ \dot
 H_p^a(\overline {\Bbb R}^n_+)$, has the unique solution
$$
u=\Xi _+^{-a}e^+r^+\Xi ^{-a}_-e^+f,\tag2.26
$$ 
lying in $
\Xi _+^{-a}(e^+\overline H_p^{t+a}({\Bbb
R}^n_+))=H_p^{a(t+2a)}(\overline{\Bbb R}^n_+)$, the $a$-transmission space.
\endproclaim

It is of course more difficult to treat variable-coefficient operators
on curved domains. For such cases, 
the following result
was shown in \cite{G15}:  

\proclaim{ Theorem 2.5}  Let $P$ be a classical strongly elliptic
  $\psi $do on ${\Bbb R}^n$ of order $2a>0$
  with even symbol (i.e., $P$ satisfies Hypothesis {\rm 1.1}),
 and let $\Omega $ be a smooth bounded subset of ${\Bbb R}^n$. Let $s>a-1/p'$. The  
homogeneous Dirichlet  problem {\rm (2.8)}, 
considered for 
$u\in \dot H_p^{a-1/p'+\varepsilon }(\comega)$, satisfies:
$$
f\in \overline
H_p^{s-2a}(\Omega )\implies u\in H_p^{a(s)}(\overline\Omega ),\text{
  the $a$-transmission space.}\tag2.27
$$
 
Moreover, the 
 mapping from $u$ to $f$ is a Fredholm mapping:
$$
 r^+  P\colon H_p^{a(s)}(\overline\Omega )\to \overline
H_p^{s-2a}(\Omega )\text{ is Fredholm}.\tag2.28
$$

A corollary for $s\to\infty $ is:
$$
r^+  P\colon \E_a(\comega)\to C^\infty (\comega)\text{ is Fredholm}.\tag2.29
$$

\endproclaim
 
The big step forward by this theorem is that
it describes the domain spaces in an exact way,
and shows that they depend only on $a, s,p$, not on the
operator $P$; and this works for all $s>a-1/p'$.

The argumentation involves a reduction to problems belonging to the calculus of
Boutet de Monvel, which is  described e.g.\ in \cite{B71,G96,S01,G09}. We use
techniques established more recently than \cite{H65,H85}, in
particular from \cite{G90}.
 The basic idea is to reduce the operator, on boundary patches, to
 the form
$$
P\sim \Lambda _-^{(a)}Q\Lambda _+^{(a)},\tag2.30
$$
where $\Lambda _\pm^{(a)}$ are order-reducing pseudodifferential operators, preserving support in
$\comega$ resp.\ ${\Bbb R}^n\setminus\Omega $, and $Q$ is of order 0
and satisfies the 0-transmission condition, hence belongs to the
Boutet de Monvel calculus.
 We shall not dwell on the proof here, but go on to some further
 developments of the theory.

\example{Remark 2.6} The assumption that the $\psi $do $P$ is even,
was made for simplicity, and could everywhere be replaced by the
assumption that $P$ has the $a$-transmission property with respect to
the chosen domain $\Omega $.  
\endexample

\example{Remark 2.7}In \cite{G14} (written after \cite{G15}), the results are extended to many other scales of spaces,
such as Besov spaces $B^s_{p,q}$ and Triebel-Lizorkin spaces $F^s_{p,q}$. Of
particular interest is the scale $B^s_{\infty ,\infty }$, also denoted
$C^s_*$, the H\"older-Zygmund scale. Here $C^s_*$ identifies with the
H\"older space $C^s$ when $s\in
\rp\setminus {\Bbb N}$, and for positive integer $k$ satisfies
$ C^{k-\varepsilon }\supset C^k_*\supset 
C^{k-1,1}\supset  C^k$ for small $\varepsilon >0$; moreover,
$C^0_*\supset L_\infty $. Then Theorem 2.5 holds with $H^s_p$-spaces replaced
by $C^s_*$-spaces.
\endexample

\example{Remark 2.8} The above applications of pseudodifferential
theory require that the domain $\Omega $ has $C^\infty $-boundary;  
in comparison,
the  results of e.g.\ Ros-Oton and coauthors in low-order H\"older spaces
allow low regularity of $\Omega $, using rather different methods.
There exists a pseudodifferential theory with just H\"older-continuous
$x$-dependence (see e.g.\ Abels \cite{A05,A05a} and references), which
may be useful to reduce the present smoothness
assumptions, but non-smooth coordinate changes for $\psi $do's 
have not yet (to our knowledge)
been established in a sufficiently useful way. At any rate, the
results obtainable by $\psi $do methods can serve as a guideline for
what one can aim for on domains with lower smoothness.  

\endexample

\head{3. Nonhomogeneous boundary value problems}\endhead

When solutions $u$ of the homogeneous Dirichlet problem lie in $d^a$
times a Sobolev or H\"older space over $\overline\Omega $, there is a
boundary value $\gamma _0(u/d^a)$, denoted 
$$
\gamma _0^au=\gamma _0(u/d^a); \tag3.1
$$
 it is
viewed as a Neumann boundary value. (We omit normalizing constants for
now; they are decribed precisely in Remark 3.2 below.)

Ros-Oton and Serra \cite{RS14a,RS15} showed the following integration-by-parts
formula:

\proclaim{ Theorem 3.1} When $u$ and $u'$ are solutions of the
  homogeneous Dirichlet problem {\rm (2.8)} for $(-\Delta )^a$ on $\Omega $ with
  $f,f'\in L_\infty (\Omega )$, $\Omega $ being $C^{1,1}$, $a>0$, 
then
$$
 \int_{\Omega }((-\Delta )^au\,\partial_j \bar u'
+\partial_ju\,(-\Delta )^a \bar u')\,dx=c\int_{\partial\Omega }\nu
_j(x)\,\gamma ^a_0u\,\gamma ^a_0\bar u'\, d\sigma ;\tag3.2
$$
here $\nu =(\nu _1,\dots,\nu _n)$
is the normal vector at $\partial\Omega  $.
\endproclaim 

It is  equivalent to a certain Pohozaev-type formula, which has important applications to uniqueness questions in nonlinear
problems for $(-\Delta )^a$. It was generalized to other related
$x$-independent singular integral operators in \cite{RSV17} (with
Valdinoci), and we extended it to $x$-dependent $\psi $do's in
\cite{G16a}. (See the survey \cite{R18} for an introduction to fractional
Pohozaev identities and their applications.)

Note that the collected order of the operators in the integral over
$\Omega $ is $2a+1$; the formula generalizes a well-known formula for
$a=1$ where $\gamma ^a_0u$ is
 replaced by the Neumann trace $\gamma _0(\partial _\nu  u)$, and the
 Dirichlet trace
 $\gamma _0u$ is $0$.

What should a nonzero Dirichlet trace be in the context of fractional Laplacians? Look at the smoothest space:  
$$
\Cal E_a(\overline{\Bbb R}^n_+)=\{u=e^+x_n^av\mid v\in C^\infty (\overline{\Bbb
  R}^n_+)\}.\tag3.3$$

By a Taylor expansion of $v$,
$$
u(x)=x_ n^av(x',0)+x_n^{a+1}\partial_nv(x',0)+\tfrac12
x_n^{a+2}\partial_n^2v(x',0)+\dots\text{ for }x_n>0.\tag3.4
$$
 
If $u\in  \Cal E_{a-1}(\crnp)$, i.e.,  $u=e^+x_n^{a-1}w$ with $w\in C^\infty ( \overline{\Bbb
  R}^n_+)$, we have analogously: 
$$
u(x)=x_ n^{a-1}w(x',0)+x_n^{a}\partial_nw(x',0)+\tfrac12
x_n^{a+1}\partial_n^2w(x',0)+\dots\text{ for }x_n>0.\tag3.5
$$
 
Here $x_n^{a-1}w(x',0)$ is the only structural difference between
(3.4) and (3.5).

This term defines the Dirichlet trace: When $u\in  \Cal E_{a-1}( \overline{\Bbb
  R}^n_+)$, the Dirichlet trace is
$$
\gamma ^{a-1}_0u=\gamma _0(u/x_n^{a-1}),\text{ equal to }\gamma _0w.\tag3.6
$$
(Again we omit a normalizing constant.) 
 
It is now natural to define a  Neumann trace on 
$ \Cal E_{a-1}( \overline{\Bbb R}^n_+)$ from the second term in (3.5),
by
$$
\gamma ^{a-1}_1u=\gamma _1(u/x_n^{a-1}),\text{ equal to }\gamma _1w=\gamma _0(\partial_nw).\tag3.7
$$
Note that it equals $\gamma _0^au$ if $u\in \E_a(\crnp)$.

\example{Remark 3.2} Also higher order traces are defined on $ \Cal E_{a-1}( \overline{\Bbb R}^n_+)$,
namely the functions $\partial_n^kw(x',0)$ in (3.5). With the
correct normalizing constants they are:
$$
\gamma _k^{a-1}u=\Gamma (a+k)\gamma
_0(\partial_n^k(u/x_n^{a-1})),\quad k\in{\Bbb N}_0.\tag3.8
$$
There are analogous definitions with $a-1$ replaced by $a-M$, $a\in
{\Bbb R}_+$ and 
$M\in{\Bbb N}_0$;  see details in  \cite{G15}, in particular Th.\ 5.1
showing mapping properties, and Th.\ 6.1 showing Fredholm solvability.
For $(-\Delta )^a$ in the case where $\Omega $ is the unit ball
in ${\Bbb R}^n$, related definitions are given by  Abatangelo, Jarohs
and Saldana in \cite{AJS18}, with explicit solution formulas.

\endexample

The above definitions can be carried over to $\Omega $ (where $x_n$ is
replaced by $d(x)$), and they extend to
$H_p^{(a-1)(s)}(\overline\Omega )$ spaces for sufficiently large $s$,
cf.\ \cite{G15}.
  
\medskip
Now consider  a $P$ satisfying Hypothesis 1.1.
We can define the {\it nonhomogeneous
  Dirichlet problem} for functions  $u \in H_p^{(a-1)(s)}(\overline\Omega
  )$ (hence supported in $\comega$), by
$$
r^+Pu=f\text{ in }\Omega ,\quad \gamma _0^{a-1}u=\varphi \text{ on }\partial\Omega .\tag3.9
$$
For this we have the solvability result (\cite{G15,G14}):

\proclaim{ Theorem 3.3}   For $s>a-1/p'$,
$$
\{r^+P, \gamma ^{a-1} _0\}\colon H_p^{(a-1)(s)} (\overline\Omega )\to
\overline H_p^{s-2a} (\Omega )\times B_p^{s-a+1/p'} (\partial\Omega ) \tag3.10
$$
is a Fredholm mapping.
\endproclaim

 Here $ B_p^{s-a+1/p'} (\partial\Omega )$ is the Besov space that
 usually appears as the range space for the standard Dirichlet trace
 operator $\gamma _0$ applied to $\overline
 H_p^{s-a+1} (\Omega )$. As in (2.20) (with $a$ replaced by $a-1$), $H_p^{(a-1)(s)}(\comega)\subset \dot H_p^s(\comega)+e^+d^{a-1}\ol
H_p^{s-a+1}(\Omega )$, when $s-a+1/p'\in \rp\setminus{\Bbb N}$.

When $a<1$, the factor $d(x)^{a-1}$ is unbounded, and the solutions of
the form $u=d^{a-1}v$, for a nice  $v$ with
nonzero boundary value, blow up at $\partial\Omega $ (a detailed
analysis is given in \cite{G14} Rem.\ 2.10). Such solutions are called
``large solutions'' in Abatangelo \cite{A15}.
 Nevertheless, $u\in L_p(\Omega )$ if $1<p<1/(1-a)$.

Nonhomogeneous Dirichlet problems (also with consecutive sets of
boundary data) are considered in
\cite{G15,G14,G16a,A15} and the recent \cite{AJS18}.

\smallskip
We can moreover consider a boundary value problem where Neumann data are prescribed:
$$
r^+Pu=f\text{ in }\Omega ,\quad \gamma _1^{a-1}u=\psi \text{ on }\partial\Omega ,\tag3.11
$$
for $u\in  H^{(a-1)(s)}(\overline\Omega )$, $s>a+1/p$. 
(The boundary condition here is {\it local}; there have also been defined other,
nonlocal Neumann problems, see the overview in \cite{G16} Sect.\ 6.)
To discuss the solvability of (3.11) we can construct a
Dirichlet-to-Neumann operator  \cite{G18a}: 

\proclaim{Theorem 3.4} Let $K_D$ be a parametrix of the mapping 
$$
z\mapsto \gamma _0^{a-1}z,\text{ when }r^+Pz=0\text{ in }\Omega ,\tag3.12
$$
(an inverse when {\rm (3.10)} is a bijection). Then the mapping
$$
S_D=\gamma _1^{a-1}K_D,\tag3.13
$$
the {\bf Dirichlet-to-Neumann operator},
is a classical pseudodifferential operator of order $1$ on
$\partial\Omega $,
with principal symbol $s_{DN,0}$ derived from the principal symbol of $P$.
\endproclaim

In particular, $s_{DN,0}(x',\xi ')$ is proportional to $|\xi '|$ for
$|\xi '|\ge 1$, when
$P=(-\Delta )^a$, considered in local coordinates at the boundary.

And then we have:

\proclaim{Theorem 3.5}
When $S_{DN}$ is elliptic (i.e., $s_{DN,0}(x',\xi ')$ is invertible
for $|\xi '|\ge 1$), the Neumann problem {\rm (3.11)} satisfies:
$$
\{r^+P, \gamma ^{a-1} _1\}\colon H_p^{(a-1)(s)} (\overline\Omega )\to
\overline H_p^{s-2a} (\Omega )\times B_p^{s-a-1/p} (\partial\Omega )\tag3.14
$$
is a Fredholm mapping, for $s>a+1/p$.

\endproclaim

Note that the ellipticity holds in the case where $P=(-\Delta )^a$.

It is remarkable  that both the Dirichlet and the Neumann boundary operators
$\gamma ^{a-1}_0$    and $\gamma ^{a-1}_1$ are {\bf local}, in spite
of the nonlocalness of the operator $P$.   
\medskip

The integration by parts formula (3.1), and the consequential Pohozaev
formulas, hold on $H_p^{a(s)}$-spaces,
where the Dirichlet trace $\gamma _0^{a-1}u$ vanishes and (consequently) the Neumann
trace $\gamma _1^{a-1}u$
identifies
with $\gamma _0^{a}u$. 

The papers \cite{RS14a,RS15,RSV17,A15,G16a} did not
establish formulas where both  $\gamma _0^{a-1}u$ and $\gamma
_1^{a-1}u$ can be nonvanishing.
However, in comparison with the standard Laplacian $\Delta $, it is natural to ask whether
there are formulas generalizing the well-known full Green's formula
with
nonzero Dirichlet and Neumann data, to these operators.
This question was answered in \cite{G18a}, where we showed:

\proclaim{Theorem 3.6} When $u,v\in  H^{(a-1)(s)}(\overline\Omega ) $, then
  for $s>a+\frac12$,
$$\aligned
\int _\Omega &(Pu\,\bar v-u\,\overline {P^*v})\,dx\\ 
&=c_a\int_{\partial\Omega
}(s_0\gamma ^{a-1}_1u\,\gamma ^{a-1}_0\bar v-s_0\gamma
^{a-1}_0u\,\gamma ^{a-1}_1\bar v+B\gamma ^{a-1}_0u\,\gamma
^{a-1}_0\bar v)\, dx';
\endaligned\tag3.15
$$
here $s_0(x')=p_0(x',\nu (x'))$ for $x'\in\partial\Omega $, and
$B$ is a first-order $\psi $do on $\partial\Omega $.
\endproclaim

Note that the only term in the right-hand side that may not be local, is
the term with $B$, nonlocal in general. A closer study (work in progress) shows that $B$
vanishes if $P=(-\Delta )^a$; we also find criteria under which $B$ is
local.

We shall not begin here to describe the method of proof; it consists
of delicate localized $\psi $do considerations using the
order-reduction operators, and elements of the Boutet de Monvel calculus.

\head{4. Heat equations}\endhead

\subhead{4.1 Anisotropic spaces, H\"older estimates, counterexamples
to high spatial\linebreak
regularity}\endsubhead

For a given lower semibounded operator
$A$ in $x$-space it is of interest to study  evolution problems with a time-parameter $t$, 
$$
Au(x,t)+\partial_tu(x,t)=f(x,t)\text{ for }t>0,\quad u(x,0)=u_0(x).\tag4.1
$$
Through many years, semigroup methods, as originally presented in Hille and Phillips \cite{HP57}, have
been developed along with other methods from functional analysis
to give
interesting results for operators $A$ acting like the Laplacian and other
elliptic differential operators. Also nonlinear questions have been
treated, e.g.\ where $f$ or $A$ are allowed to depend on $u$. 

For $A$ representing the fractional
Laplacian and its generalizations, on ${\Bbb R}^n$ or a domain, the studies have begun more
recently. A natural approach is here: To find the appropriate general strategies
from the works on differential operators,
and show the appropriate properties of $A$
assuring that the methods can be applied to it.

Consider the evolution problem (heat equation) associated with an
operator $P$ as studied in Sections 2 and 3, with a homogeneous Dirichlet condition:
$$\aligned
Pu+\partial_tu&=f\text{ on }\Omega \times I ,\quad I=\,]0,T[\,,\\
u&=0\text{ on }(\Bbb R^n\setminus\Omega )\times I, \\
u|_{t=0}&=u_0.
\endaligned \tag4.2
$$

Since $P_{\operatorname{Dir}}$ is a positive
selfadjoint (or sectorial) operator in $L_2(\Omega )$, there is  
solvability in a framework of $L_2$-Sobolev spaces.

We are interested in the {\it regularity} of solutions.

This question has been treated recently by Leonori, Peral, Primo and
Soria \cite{LPPS15} in $L_r(I;L_q(\Omega ))$-spaces, by Fernandez-Real
and Ros-Oton \cite{FR17} in anisotropic H\"older spaces, and by
Biccari, Warma and Zuazua \cite{BWZ18} for $(-\Delta )^a$ 
in local $L_p$-Sobolev spaces over
$\Omega $. Earlier results are shown e.g.\ in Felsinger and Kassmann
\cite{FK13} and Chang-Lara and Davila
\cite{CD14} (H\"older properties), and Jin and Xiong \cite{JX15}
(Schauder estimates). The references in the mentioned works give
further information, also on related {\it  heat kernel} estimates.
Very recently (November 2017), the H\"older estimates were improved by
Ros-Oton and Vivas \cite{RV18}.

We have a few contributions to this subject, that we shall describe
in the following.

Let us first introduce {\bf anisotropic spaces} of Sobolev or H\"older
type. Let $d\in{\Bbb R}_+$. There are the Bessel-potential types, for $s\in{\Bbb R}$:
$$\aligned
H_p^{(s,s/d)}(\Bbb R^n\!\times\!\Bbb R)&=\{u\in \Cal
S'\mid \Cal F^{-1}((\langle\xi \rangle^{2d}+\tau
^2)^{s/2d}\hat u(\xi ,\tau ))\in L_p(\Bbb R^{n+1})\},\\ 
\overline H_p^{(s,s/d)}(\Omega \times I)&=r_{\Omega \times
  I}H_p^{(s,s/d)}(\Bbb R^n\!\times\!\Bbb R),
\endaligned\tag4.3$$
and there are related definitions of Besov-type with $H_p$ replaced by
$B_p$.
There are the  H\"older spaces:
$$
\overline C^{(s,r)}(\Omega \times I)=L_\infty
(I;\overline C^{s}(\Omega ))\cap L_\infty (\Omega ;\overline C^{r}( I)),
\text{ for }s,r\in \Bbb R_+;\tag4.4
$$
including in particular the case $r=s/d$. (For $s$ equal to an integer
$k$, we denote by $\ol C^k(\Omega )$ the space of bounded continuous
functions on $\comega$ with bounded derivatives up to order $k$, this
includes the case $\Omega ={\Bbb R}^n$.)
The spaces occur in many works; important properties
are recalled e.g.\ in \cite{G95} and \cite{G18}
with further references. 

The H\"older-type spaces are the primary
objects in the investigations of
Fernandez-Real and Ros-Oton in \cite{FR17}, Ros-Oton and Vivas in \cite{RV18}.
These authors have for the Dirichlet heat problem  the following
results, showing the role of $d^a$ in H\"older estimates:

\proclaim{Theorem 4.1}
Let $P$ be an operator of
the form {\rm (1.2)ff.}, $0<a<1$, and consider solutions of the problem {\rm (4.2)}.

$1^\circ$ {\rm \cite{FR17}, Cor.\ 1.6.} When $\Omega $ is a bounded open $C^{1,1}$ subset
of ${\Bbb R}^n$, then the unique weak solution $u$ 
 with
$f\in L_\infty (\Omega \times I)$ and $u_0\in L_2(\Omega )$ satisfies
$$
\aligned
\|u\|_{\ol C^{(a,1-\varepsilon )}(\Omega \times I')}&+\|u/{d^a}\|_{\ol
C^{(a-\varepsilon ,1/2-\varepsilon/(2a)) }(\Omega \times I')}\\
&\le C
(\|f\|_{L_\infty (\Omega \times I)}+\|u_0\|_{L_2(\Omega )}),
\endaligned\tag4.5
$$
for any small $\varepsilon >0$, $I'=\,]t_0,T[\,$ with $t_0>0$.
Moreover, if $f\in \ol C^{(\gamma,\gamma/(2a))}(\Omega \times I)$ with $\gamma\in\,]0,a]$
such that $\gamma+2a\notin{\Bbb N}$, then $u$ has the interior regularity:
$$
\|u\|_{\ol C^{(2a+\gamma,1+\gamma/(2a)) }(\Omega '\times I')}\le C'\|f\|_{ C^{(\gamma,\gamma/(2a))}(\Omega \times I)},\tag4.6
$$
for any
 $\Omega '$ with $\ol {\Omega '}\subset \Omega $.

$2^\circ$ {\rm \cite{RV18}, Cor.\ 1.2.} Let $\gamma \in \,]0,a[\,$, $\gamma +a\notin{\Bbb N}$. 
When $\Omega $ is a bounded open $C^{2,\gamma  }$ subset
of ${\Bbb R}^n$, then 
$f\in \ol C^{(\gamma ,\gamma /(2a))}(\Omega \times I)$, $u_0\in L_2(\Omega )$ imply:
$$
\aligned
\|u\|_{\ol C^{(\gamma ,1+\gamma /(2a) )}(\Omega \times I')}&+\|u/{d^a}\|_{\ol
C^{(a+\gamma  ,1/2+   \gamma /(2a)) }(\Omega \times I')}\\
&\le C
(\|f\|_{\ol C^{(\gamma ,\gamma /(2a))} (\Omega \times
I)}+\|u_0\|_{L_2(\Omega )}),
\endaligned\tag4.7
$$
for any $I'=\,]t_0,T[\,$ with $t_0>0$.

\endproclaim

In other words, if we take $\gamma =a-\varepsilon $ for a small
$\varepsilon >0$, (4.7) reads, with $\varepsilon '=\varepsilon /(2a)$,
$$
\aligned
\|u\|_{\ol C^{(a-\varepsilon  ,3/2-\varepsilon ' )}(\Omega \times I')}&+\|u/{d^a}\|_{\ol
C^{(2a-\varepsilon  ,1-\varepsilon ') }(\Omega \times I')}\\
\le C
&(\|f\|_{\ol C^{(a -\varepsilon ,1/2-\varepsilon ')} (\Omega \times
I)}+\|u_0\|_{L_2(\Omega )}). \endaligned \tag4.8
$$

An interesting question is whether the regularity of
$u$  can be lifted further, when $f$ in $2^\circ$ is replaced by a more
regular function. As we shall see below, this is certainly possible
with respect to the $t$-variable; this is also shown to some extent in
\cite{FR17}. However, there are limitations with respect to the
boundary behavior in the
$x$-variable.
It is shown in \cite{G15a} that when $P$ is as in Hypothesis 1.1, any eigenfunction $\varphi $ of
$P_{\operatorname{Dir}}$  associated with a nonzero eigenvalue $\lambda
$ satisfies (more on the spaces in Theorem 4.2):$$
\varphi \in C_*^{a(3a)}(\comega )\subset d^a\ol C_*^{2a}(\comega)\; \cases = d^a\ol C^{2a}(\Omega )\text{ if
}a\ne \frac12,\\ \subset d^a\ol C^{2a-\varepsilon }(\Omega )\text{  if }
a=\frac12,\endcases\tag4.9$$
but, in the basic case  $P=(-\Delta )^a$,  $\varphi $ is {\it not in}  
$d^a\ol C^\infty (\Omega )$ and not either in $\ol C^\infty (\Omega
)$. The function $u(x,t)=e^{-\lambda
t}\varphi (x)$ is clearly a solution of the heat equation with $f=0$ (hence $f\in \ol C^\infty (I;\ol C^\infty (\Omega ))$),
but $u$ and $u/d^a\notin L_\infty (I';\ol C^\infty (\Omega ))$. This
shows a surprising contrast to the usual 
regularity  rules for heat
equations, and it differs radically from the stationary case, where we have (2.29).

The argument can be extended from $(-\Delta )^a$ to more general operators, and it can be sharpened to rule out also
finite higher order regularities.

\proclaim{Theorem 4.2} Let $P$ satisfy Hypothesis {\rm 1.2} with $0<a<1$,
or let $P$ equal $(-\Delta )^a$ with $a>0$, or the fractional Helmholtz
operator $(-\Delta +m^2)^a$ with $m>0$ and $0<a<1$. Let $\Omega $ be $C^\infty $. Then
any eigenfunction $\varphi $ of $P_{\operatorname{Dir}}$ associated
with a nonzero eigenvalue $\lambda $ satisfies {\rm (4.9)}, but is not in
$\overline C^{a+\delta }(\Omega )$ nor in $C_*^{a(3a+\delta )}(\comega)$
for any $\delta >0$.
\endproclaim

\demo{Proof} Recall  that $C^s_*$ stands for the
H\"older-Zygmund space, which identifies with $C^s$ when $s\in
\rp\setminus {\Bbb N}$, cf.\ Remark 2.7.
As shown in
\cite{G14}, p.\ 1655 and Th.\ 3.2,  $C_*^{a(2a+s )}(\comega)$ is the solution
space for the homogeneous Dirichlet problem with right-hand side in
$\ol C_*^{s }(\Omega )$; here $C_*^{a(2a+s
)}(\comega)\subset d^a\ol C_*^{a+s }(\Omega )$, but there is not
equality.
(One reason for the lack of equality is that the functions in $C_*^{a(2a+s )}(\comega)$ are in
$C_*^{2a+s }$ over the interior, another is that $C_*^{a(2a+s
)}(\comega)$ only reaches a subspace
of $d^a\ol C_*^{a+s }(\Omega )$ near the boundary, cf.\ \cite{G15} Th.\ 5.4.) 

Assume that  $\varphi $ satisfies  $P_{\operatorname{Dir}}\varphi
=\lambda \varphi $ ($\lambda \ne 0$) and is in $\ol C^{a+\delta }(\Omega )$ for a positive
$\delta <1-a$; then in fact $\varphi \in \dot C^{a+\delta }(\comega)$
since $\varphi \in C_*^{a(3a+\delta )}(\comega)\subset d^a\ol
C_*^{2a+\delta }(\Omega )$ implies $\gamma _0\varphi =0$. Now
$\varphi \in \dot C^{a+\delta }(\comega )$ implies $\gamma _0(\varphi
/d^a)=0$ since $\delta >0$.
 It is shown in \cite{RSV17} for operators of the
form (1.2)ff.\ with $0<a<1$, in \cite{RS15} for $(-\Delta )^a$ with $a>0$, and in \cite{G16a}, Ex.\ 4.10  for $(-\Delta +m^2)^a$,
how it follows from Pohozaev identities that
$$
P_{\operatorname{Dir}}v=\lambda v, \; \gamma _0(v/d^a)=0\implies
v\equiv 0.
$$
 Thus $\varphi =0$ and cannot be an eigenfunction.\qed
\enddemo

This allows us to conclude:

\proclaim{Corollary 4.3}
Consider the problem {\rm (4.2)}. 
For the operators $P$ considered in Theorem {\rm 4.2}, $C^\infty
$-regularity of $f$ does not imply $C^\infty $-regularity of $u$ or $u/d^a$. In
fact, there exist
choices of $f(x,t)\in \ol C^\infty (\Omega
\times I )$ with
solutions $u(x,t)$ satisfying
$$
u\notin L_\infty (I'; \overline C^\infty (\Omega )),\quad u\notin
L_\infty (I';d^a \ol C^\infty (\Omega)).\tag4.10
$$

More precisely, there exist solutions with  $f(x,t)\in \ol C^\infty (\Omega
\times I )$  such that
$$
u\notin L_\infty (I'; \overline C^{a+\delta }(\Omega )),\quad u\notin
L_\infty (I';C_*^{a(3a+\delta )}(\comega)), \text{ for any $\delta >0$.  }\tag4.11
$$

\endproclaim

\demo{Proof} This follows by taking  $u(x,t)=e^{-\lambda t}\varphi (x)$
with an eigenfunction $\varphi $ as in Theorem 4.2. It solves the heat
problem (4.2) with $f=0$ (thus $f\in
 \ol C^\infty (\Omega \times
 I)$) and $u_0=\varphi $, and it clearly satisfies (4.10) as well as (4.11).\qed
\enddemo

Note that the $x$-regularity obtained in (4.8) is close to the
upper bound.

Corollary 4.3 gives counterexamples; one can show more systematically
that $\gamma ^a_0u$ being nonzero prevents the solution from being in
$C^\infty $ or $d^aC^\infty  $ at the boundary \cite{G18b}.

\subhead 4.2 Solvability in Sobolev spaces \endsubhead

Now let us consider estimates in Sobolev spaces.

In \cite{GS90} (jointly with Solonnikov) and in \cite{G95} the author studied
evolution problems for $\psi $do's $P$ with the 0-transmision property
at $\partial\Omega $, along with trace, Poisson and singular Green
operators in the Boutet de Monvel calculus (cf.\ \cite{B71,G90,G96}),
setting up a full calculus leading to existence, uniquenes and
regularity theorems in anisotropic Bessel-potential and Besov spaces as
mentioned 
in (4.3)ff.

These works take $P$ of integer order, and do not cover the present
case. We expect that a satisfactory generalization of the full boundary
value theory in those works, to heat problems for
our present operators, would be quite difficult to achieve.
However their point of view on the $\psi $do $P$ alone, considered on
${\Bbb R}^n$ without boundary conditions, can be extended, as follows:

For a  classical
strongly elliptic $\psi $do $P$ of order $d\in\rp$ on $\Bbb R^n$ (with global
symbol e\-sti\-ma\-tes), we can construct an anisotropic symbol
calculus on ${\Bbb R}^n\stimes{\Bbb R}$ that includes operators
$P+\partial_t$ and their parametrices. It is not
quite standard, since the typical strictly  homogeneous symbol $|\xi |^d+i\tau $ is not
$C^\infty $ at points $(0,\tau )$ with $\tau \ne 0$. But this is a
phenomenon handled in \cite{G96} by introducing classes of symbols
with finite ``regularity number'' $\nu$ (essentially the H\"older
regularity of the strictly homogeneous principal symbol at points
$(0,\tau )$), and keeping track of how  the value of $\nu $ behaves in compositions
and parametrix constructions.
 
The calculus gives, on $\Bbb R^{n+1}$ and locally in $\Omega \times I$ \cite{G18}:

\proclaim{Theorem 4.4}  Let $P$ be a classical strongly elliptic
  $\psi $do of order $d\in{\Bbb R}_+$. Then $P+\partial_t$ maps $H^{(s,s/d)}_p(\Bbb
  R^n\!\times\!\Bbb R)$ continuously into $H^{(s-d,s/d-1)}_p(\Bbb
  R^n\!\times\!\Bbb R)$ for any $s\in {\Bbb R}$. Moreover:

$1^\circ$ If $u\in
H^{(r,r/d)}_p(\Bbb R^n\!\times\!\Bbb R) $ for some large negative $r$
(this holds in
particular if $u\in\E'({\Bbb R}^{n+1})$ or e.g.\ $L_p({\Bbb
R};\E'({\Bbb R}^n))$),  then  
$$
(P+\partial_t)u\in H^{(s,s/d)}_p(\Bbb R^n\!\times\!\Bbb
R)\implies
u\in H^{(s+d,s/d+1)}_p(\Bbb R^n\!\times\!\Bbb R).\tag4.12
$$

$2^\circ$ Let $\Sigma =\Omega \times I$, and let
$u\in
H^{(s,s/d)}_p(\Bbb R^n\!\times\!\Bbb R) $. Then
$$
(P+\partial_t)u|_{\Sigma }\in H^{(s,s/d)}_{p,\operatorname{loc}}(\Sigma )\implies
u\in H^{(s+d,s/d+1)}_{p,\operatorname{loc}}(\Sigma ).\tag4.13
$$ 
\endproclaim 

This theorem works for {\it any} strongly elliptic  classical $\psi
 $do of positive order, not just fractional Laplacians, but for example also $-\Delta
+(-\Delta )^{1/2}$ or $(-\Delta )^{1/2}+b(x)\cdot\nabla +c(x)$ with
 real $C^\infty $-coefficients. The result extends by standard
 localization methods to the case where ${\Bbb R}^n$ is replaced by a closed manifold.

Note that in $2^\circ$, the regularity of $u$ is only lifted by 1 in
$t$, and the hypothesis on $u$ concerns all $x\in \Bbb R^n$; the
necessity of this is pointed out in related situations in \cite{CD14}
and \cite{FR17}.

By use of embedding theorems, we can moreover derive from the above a local regularity
result in anisotropic H\"older spaces:

\proclaim{Theorem 4.5} Let $P$ and $\Sigma $ be as in Theorem {\rm
4.3.}  Let $s\in \rp$, and let $u\in 
C^{(s,s/d)}(\Rn\!\times\!\Bbb R)\cap \E'(\Rn\!\times\!\Bbb R)$.
Then
$$
(P+\partial_t)u|_{\Sigma  }\in
C^{(s,s/d)}_{\operatorname{loc}}(\Sigma )\implies u|_{\Sigma }\in
 C^{(s+d-\varepsilon ,(s-\varepsilon )/d+1)}_{\operatorname{loc}}(\Sigma ),\tag 4.14
$$
for small $\varepsilon >0$.
\endproclaim

Observe the similarity with (4.6) (where $d=2a$). On one hand we have a loss of
$\varepsilon $; on the other hand we have general $x$-dependent
operators $P$ just required to be strongly elliptic (albeit with
smooth symbols), and no upper limitations on $s$. The $\varepsilon $ might
possibly be removed by working with our operators on the H\"older-Zygmund
scale
(cf.\ e.g.\ \cite{G14}).

Still other spaces could be examined. There is the work of Yamazaki \cite{Y86} on
$\psi $do's acting in  anisotropic Besov-Triebel-Lizorkin spaces
(defined in his work);
these  include the $H_p$ and $B_p$
spaces as special cases. However, the operators in \cite{Y86} seem to be more
regular 
(their quasi-homogeneous symbols being smooth outside of 0) than the operators $\partial_t+P$ that we
study here.
There is yet another type of spaces, the so-called modulation spaces on ${\Bbb R}^n$, where there
are very recent results on heat equations (with nonlinear
generalizations)  by Chen, Wang, Wang and Wong \cite{CWWW18}.

For the case where  boundary conditions at $\partial\Omega $ are imposed,
there is not (yet) a systematic boundary-$\psi $do theory as in the
stationary case. But using suitable functional analysis results we can
make some progress, showing how the heat equation solutions behave in  terms
of  Sobolev-type spaces involving the factor $d^a$.
 We henceforth restrict the
attention to the case $0<a<1$, and to problems with initial value 0,
$$\aligned
Pu+\partial_tu&=f\text{ on }\Omega \times I ,\quad I=\,]0,T[\,,\\
u&=0\text{ on }(\Bbb R^n\setminus\Omega )\times I, \\
u|_{t=0}&=0.
\endaligned \tag4.15
$$

There is a straightforward result in the $L_2$-framework:

\proclaim{Theorem 4.6}  Let $P$ satisfy Hypothesis {\rm 1.1} with $a<1$,
and let $\Omega $ be a smooth bounded subset of
  ${\Bbb R}^n$. For $f$ given in $L_2(\Omega \times I)$, there is a
  unique solution $u$ of {\rm (4.15)} satisfying
$$
u\in L_2(I;H^{a(2a)}(\comega ))\cap \ol H^1(I; L_2(\Omega ));\tag4.16
$$
 here $H^{a(2a)}(\comega )=D(P_{\operatorname{Dir},2})$ equals $\dot
 H^{2a}(\comega)$ if $a<\frac12$, and is as described in {\rm
 (2.13)}ff.\  (with $p=2$) when $a\ge \frac12$. 

Moreover, $u\in \ol C^0(I;L_2(\Omega ))$.

\endproclaim

\demo{Proof} 
Define the sesquilinear form
$$
Q_0(u,v)=(r^+Pu,v)_{L_2(\Omega )},
$$
 first for $u,v\in C_0^\infty (\Omega )$ and then extended by closure in
$\dot H^a$-norm to a bounded sesquiliear form with domain $\dot
H^a(\comega)$. In view of the strong ellipticity, it is coercive:
$$
\operatorname{Re}Q_0(u,u)\ge c_0\|u\|^2_{\dot H^a}-\xi \|u\|^2_{L_2}
\text{ with }c_0>0,\xi \in{\Bbb R},\tag4.17
$$
and hence defines via the Lax-Milgram lemma a realization of $r^+P$
with domain (2.9); for precision we shall denote the operator  
$P_{\operatorname{Dir},2}$. (The Lax-Milgram construction is described
e.g.\ in \cite{G09}, Sect.\ 12.4.) The adjoint is defined similarly from
$Q_0^*(u,v)=\overline {Q_0(v,u)}$, and it follows from (4.17) that the
spectrum and numerical range is contained in a set 
$\{z\in{\Bbb C}\mid
|\operatorname{Im}z|\le C(\operatorname{Re}z+\xi ), \, \operatorname{Re}z>\xi _0\}$, and 
$$
\|(P_{\operatorname{Dir},2} -\lambda )^{-1}\|_{\Cal L(L_2)}\le
c\ang{\lambda }^{-1}\text{ for }\operatorname{Re}\lambda \le -\xi _0.\tag4.18 
$$
We can then apply Lions and Magenes \cite{LM68} Th.\ 4.3.2, which
shows that there is a unique solution $u$ of (4.15) in
$L_2(I;D(P_{\operatorname{Dir},2}))$. Here, moreover,
$\partial_tu=f-r^+Pu\in L_2(\Omega \times I)$, so $u\in \ol H^1(I;L_2(\Omega ))$.

The last statement follows since $ \ol H^1(I;L_2(\Omega ))\subset \ol C^0(I;L_2(\Omega ))$.\qed

\enddemo

The domain of $P_{\operatorname{Dir},2}$ is  a  Sobolev
space $\dot H^{2a}(\comega)$ when $a<\frac12$, but not a standard
Sobolev space when $a\ge \frac12$. However, it is then contained in
one. Namely, if we take $r\ge 0$ such that $$
r\;\cases =2a\text{ if }0<a<\tfrac12,\\ <a+\tfrac12 \text{ if }\tfrac12\le a<1 ,\endcases\tag4.19
$$
 then 
$$
H^{a(2a)}(\comega)\subset \dot H^r(\comega)\subset \overline
H^{r}(\Omega ),\text{ hence }\|u\|_{\ol H^r}\le c\|u\|_{H^{a(2a)}},\tag 4.20
$$
for all $0<a<1$. This follows for $\Omega =\rnp$ since $\ol
H^a(\rnp)=\dot H^a(\crnp) $ if $a<\frac12$, and  $\ol
H^a(\rnp)\subset \dot H^{\frac12-\varepsilon }(\crnp) $ if
$a\ge\frac12$, so that by  (2.19),
$$
H^{a(2a)} (\overline{\Bbb R}^n_+)=\Xi _+^{-a}e^+\overline
H^{a} ({\Bbb R}^n_+)\subset \Xi _+^{-a}\dot
H^{r-a} (\crnp)=\dot H^{r} (\crnp);\tag4.21
$$
 there is a similar proof for general $\Omega $. Observe that $r\le
 2a$ in all cases.

This can be used to jack up the regularity
result of Theorem 4.6 by one derivative in $t$ and an improved
$x$-regularity, when $f$ is $H^r$ in $x$ and $H^1$ in $t$. Higher
$t$-regularity can also be obtained. To do this, we shall apply the more refined
Th.\ 4.5.2 in \cite{LM68}, introduced there for the purpose of showing
higher regularities. For the convenience of the reader, we list a
slightly reformulated version:

\proclaim{Theorem 4.7}{\rm (From \cite{LM68} Th.\ 4.5.2.)} Let $X$ and $\Cal H$ be Hilbert
spaces, with $X\subset \Cal H$, continuous injection. Let $A$ be an unbounded
linear operator in $X$ such that $A-\lambda$ is a bijection from the domain
$
D_X(A)=\{u\in X\mid Au\in X\}
$
onto $X$, for all $\lambda\in{\Bbb C}$ with
 $\operatorname{Re}\lambda\le -\xi _0$.

Assume moreover that for all such $\lambda$, and for $u\in D_X(A)$,
$$
\|(A-\lambda)u\|_{X}+\ang \lambda^\beta \|(A-\lambda)u\|_{\Cal H}\ge c(\|u\|_{D_X(A)}+
\ang \lambda^{\beta +1}\|u\|_{\Cal H}),\tag4.22
$$
where $\beta >0$  and $c>0$ are given.

Let $\beta \notin \frac12+{\Bbb N}$. The problem $
Au+\partial_tu=f\text{ for }t\in I$, $ u(0)=0$,
with $f$ given in $L_2(I;X)\cap H^\beta (I;\Cal H)$ with $f^{(j)}(0)=0$
for $j<\beta -\frac12$, has a unique solution $$
u\in L_2(I;D_X(A))\cap H^{\beta +1}(I;\Cal H).\tag4.23
$$

\endproclaim

This allows us to show:

\proclaim{Theorem 4.8} Assumptions as in Theorem {\rm 4.6}.

$1^\circ$ If $f\in L_2(I;\ol H^{r}(\Omega ))\cap \ol H^1(I;L_2(\Omega ))$ for
some $r$ satisfying {\rm (4.19)}, with $f|_{t=0}=0$, then the
solution of {\rm (4.15)} satisfies
$$
u\in L_2(I; H^{a(2a+r)}(\comega ))\cap \ol H^2(I;L_2(\Omega )).\tag4.24
$$

$2^\circ$ For any integer $k\ge 2$, if
$
f\in L_2(I;\ol H^{r}(\Omega
))\cap \ol H^k(I;L_2(\Omega ))$ with $\partial_t^jf|_{t=0}=0$ for
$j<k$, then
$$
u \in L_2(I; H^{a(2a+r)}(\comega ))\cap \ol H^{k+1}(I;L_2(\Omega )).\tag4.25
$$
It follows in particular that
$$
 f\in \bigcap_k \ol H^k(I;\ol H^r(\Omega
 )),\;\partial_t^jf|_{t=0}=0\text{ for } j\in{\Bbb N}_0\implies
u \in \bigcap_k\ol H^{k}(I; H^{a(2a+r)}(\comega )).\tag4.26
$$
\endproclaim

\demo{Proof}
With $A$ acting
like $r^+P$, denote 
$$
D_r(A)=\{v \in H^{a(2a)}(\comega)\mid Av\in \ol H^{r}(\Omega)\}.\tag4.27
$$
Here $D_r(A)=H^{a(2a+r)}(\comega)$ in view of Theorem 2.5 (since $r\ge
0$). Moreover, it equals $D_{\ol H^r}(A)=\{v \in \ol H^r(\Omega)\mid Av\in
\ol H^{r}(\Omega)\}$, since $H^{a(2a)}(\comega)\subset \dot
H^r(\comega)\subset \ol H^r(\Omega )$. When $\operatorname{Re}\lambda
\le -\xi _0$, the bijectiveness of $A-\lambda $ from $H^{a(2a)}(\comega)$
to $L_2(\Omega )$ implies bijectiveness from $D_r(A)$ to $\ol
H^r(\Omega )$. All this shows that $D_r(A)=D_{\ol H^r}(A)$ is as in
the start of Theorem 4.7 with $X=\ol H^r(\Omega )$, $\Cal H=L_2(\Omega )$.
Moreover, there is an equivalence of norms
$$
\|(A+\xi _0)v\|_{ \ol H^{r}}\simeq\|v\|_{D_r(A)},   \text{ for }v\in D_r(A).\tag4.28
$$
Note also that besides the inequality (4.18), that may be written
$$
\ang\lambda \|(A-\lambda )^{-1}g\|_{L_2}\le
c\|g\|_{L_2}, \text{ for }g\in L_2(\Omega ),\quad
\operatorname{Re}\lambda \le -\xi _0, \tag4.29
$$
we have that by (4.20), (4.28) for $r=0$ and (4.29),
$$
\aligned
 \|(A-\lambda )^{-1}g\|_{\ol H^{r}}&\le c\|(A-\lambda
 )^{-1}g\|_{H^{a(2a)}}\simeq \|(A+\xi _0)(A-\lambda )^{-1}g\|_{L_2}\\
&\le \|g\|_{L_2}+|\lambda +\xi _0| \|(A-\lambda
 )^{-1}g\|_{L_2}\le c'\|g\|_{L_2}, \text{ for }g\in L_2(\Omega ),
\endaligned
$$
so altogether,
$$
\ang\lambda \|(A-\lambda )^{-1}g\|_{L_2}+ \|(A-\lambda )^{-1}g\|_{\ol H^{r}}\le
c_1\|g\|_{L_2}, \text{ for }g\in L_2(\Omega ).\tag4.30
$$
For $v\in D_r(A)$, with $(A-\lambda )v$ denoted $g$, (4.28) implies
$$
\aligned
\|v\|_{D_r(A)}&\simeq\|(A+\xi _0)v\|_{\ol H^r}\le c_2(\|(A-\lambda
)v\|_{\ol H^{r}} +|\lambda +\xi _0
|\|v\|_{\ol H^{r}})\\
&\le c_3(\|(A-\lambda )v\|_{\ol H^{r}}+\ang\lambda \|g\|_{L_2}),
\endaligned \tag4.31
$$
where we used (4.30) in the last step. 
Moreover, by (4.30),
$$
\ang\lambda ^2\|v\|_{L_2}\le c_1\ang\lambda \|g\|_{L_2},
$$
so we altogether find the inequality  for $v\in D_r(A)$:
$$
\|v\|_{D_r(A)}+\ang\lambda ^2\|v\|_{L_2}\le c_4(\|(A-\lambda )v\|_{\ol H^{r}}+\ang\lambda
\|(A-\lambda )v\|_{L_2}),\quad \operatorname{Re}\lambda \le -\xi _0.\tag4.32
$$

We can now apply Theorem 4.7,
with $\beta =1$, $X=\ol H^{r}(\Omega)$ and $\Cal
H=L_2(\Omega )$.
It follows that the solution
of (4.15) with $f(x,t)$ given in $L_2(I;\ol H^r(\Omega ))\cap \ol
H^1(I;L_2(\Omega ))$, $f|_{t=0}=0$, satisfies
$$
u\in L_2(I;D_r(A))\cap \ol H^2(I; L_2(\Omega )),
$$
from which (4.19) follows since $D_r(A)=H^{a(2a+r)}(\comega)$. This
shows $1^\circ$.

Now let $k\ge 2$. By (4.29) and (4.31), since $v=(A-\lambda )^{-1}g$,
$$\aligned
\|v\|_{D_r(A)}+\ang\lambda ^{k+1}\|v\|_{L_2}&\le c_3(\|(A-\lambda
)v\|_{\ol H^{r}}+\ang\lambda \|g\|_{L_2})+c\ang\lambda
^{k}\|g\|_{L_2}\\
&\le c_5(\|(A-\lambda )v\|_{\ol H^{r}}+\ang\lambda ^{k}\|(A-\lambda )v\|_{L_2}),
\endaligned\tag 4.33$$
which allows an application of Theorem 4.7 with $\beta
=k$,  $X=\ol H^{r}(\Omega)$, $\Cal
H=L_2(\Omega )$, giving the conclusion (4.25). (4.26) follows when
$k\to\infty $. This shows $2^\circ$. 
\qed

\enddemo

Note in particular that  $D_r(A)= H^{a(4a)}(\comega ))$ when
$a<\frac12$.

From the point of view of anisotropic Sobolev spaces, the solutions
in (4.24) and (4.25) satisfy, since $H^{a(2a+r)}(\comega)\subset
H^{a(2a)}(\comega)\subset \dot H^{r}(\comega)\subset \ol H^r(\Omega ) $,
$$
u\in L_2(I;\ol H^r(\Omega ))\cap \ol H^2(I;L_2(\Omega ))=\ol
H^{(r,2)}( \Omega\times I ), \text{ resp.\ }u\in \ol
H^{(r,k+1)}( \Omega\times I ),
$$
but (4.24)--(4.26) give a more refined information. 

Note that $r<3/2$ in all these cases. 
 We think that a lifting to higher values of $r$, 
 of the conclusion of $1^\circ$ concerning regularity in $x$ at the boundary, 
 would demand very different methods, or may not even be
 possible, because of the incompatibility of $H^{a(s)}(\comega)$ with standard
 high-order Sobolev spaces when $s$ is high. Cf.\ also Corollary 4.3,
 which excludes higher smoothness in a related situation.

\medskip
Next, we turn to $L_p$-related Sobolev spaces with
general $p$. The following
result was shown for 
$x$-independent operators in \cite{G18}:

\proclaim{Theorem 4.9}  Let $P$ satisfy Hypothesis {\rm 1.2} with $a<1$.
Then {\rm (4.15)}
has for any  $f\in L_p(\Omega \times I)$  a unique solution $u(x,t)\in
\ol C^0( I;L_p(\Omega ))$ ($1<p<\infty $); it  satisfies:
$$
u\in L_p(I;H_p^{a(2a)}(\overline\Omega))\cap \ol H^1_p( I;L_p(\Omega
)).\tag4.34
$$

\endproclaim
 
The result is sharp, since it gives the exact domain for $u$ 
mapped into $f\in L_p(\Omega \times I)$.

The proof relies on a theorem of Lamberton \cite{L87} which was used earlier in the work of Biccari, Warma and Zuazua
\cite{BWZ18}; they showed a local version of Theorem 4.9 for $P=(-\Delta )^a$, where
$H_p^{a(2a)}(\overline\Omega)$  in (4.34) is replaced
by $B^{2a}_{p,\operatorname{loc}}(\Omega )$ if $p\ge 2$,
$a\ne\frac12$, by $H^1_{p,\operatorname{loc}}(\Omega )$ if
$a=\frac12$, and by $B^{2a}_{p,2,\operatorname{loc}}(\Omega )$ if
$p<2$, $a\ne \frac12$.

\example{Remark 4.10} Another paper \cite{BWZ17}, preparatory for
\cite{BWZ18}, is devoted to a computational proof of local
regularity (regularity in compact subsets of $\Omega $) of the solutions of the
stationary Dirichlet problem (2.8) for $P=(-\Delta )^a$ with $f\in L_p(\Omega
)$. Here the 
authors were apparently unfamiliar with the pseudodifferential elliptic
regularity theory  (mentioned in \cite{G15}) that
gave the answer many years ago, see the addendum \cite{BWZ17a}. The
addendum also corrects some mistakes connected with the definition of
$W^{s,p}$-spaces, cf.\ (2.3)ff.\ above. The  proof of 
local regularity is repeated in \cite{BWZ18}.
\endexample

 Let us recall the proof of Theorem 4.9 from \cite{G18}: The Dirichlet realization in $L_p(\Omega
 )$, namely the operator $P_{\operatorname{Dir},p}$, acting like
 $r^+P$ with domain
$$
D(P_{\operatorname{Dir},p})=\{u\in \dot H_p^a(\overline\Omega )\mid r^+Pu\in L_p(\Omega )\},\tag4.35
$$
coincides on $L_2(\Omega )\cap
 L_p(\Omega )$ with $P_{\operatorname{Dir},2}$  defined variationally
 from the sesquilinear form
$$
Q_0(u,v)=\int_{\Bbb R^{2n}}(u(x)-u(y))(\bar v(x)-\bar
v(y))K(x-y)\,dxdy \text{ on }\dot H^a(\overline\Omega ),
$$
where $K(y)=c\Cal F^{-1}p(\xi )$, positive and homogeneous of degree $-2a-n$. 
 
The form has the Markovian property: When $u_0$ is defined from a real
function $u$ by $
u_0=\min\{\max\{u,0\},1\}$, then  $Q_0(u_0,u_0)\le
Q_0(u,u)$. It is a so-called Dirichlet form, as explained in  Fukushima, Oshima and Takeda
\cite{FOT94}, pages 4--5 and  Example 1.2.1, and Davies \cite{D89}.
 Then, by \cite{FOT94} Th.\ 1.4.1 and \cite{D89} Th.\ 1.4.1--1.4.2, $-P_{\operatorname{Dir},p}$
generates a strongly continuous contraction semigroup $T_p(t)$ not only in $L_2(\Omega )$ for $p=2$ but also
in $L_p(\Omega )$ for any $1<p<\infty $, and $T_p(t)$ is bounded holomorphic.
Then
the statements in Theorem 4.9 follow from Lamberton \cite{L87} Th.\ 1.\qed

\example{Remark 4.11} An advantage of the above results is that we
have a precise characterization of the Dirichlet domain, namely 
$D(P_{\operatorname{Dir},p})=H_p^{a(2a)}(\overline\Omega )$. This space can be
further described, as 
shown in \cite{G15}, cf.\ Th.\ 5.4.
When $p<1/a$ then $H_p^{a(2a)}(\overline\Omega )=\dot
H_p^{2a}(\overline\Omega )$. Since  $\dot
H_p^{2a}(\overline\Omega )\subset \ol H_p^{2a}(\Omega )$, the solutions $u$ are in the
anisotropic space
$$
\ol H_p^{(2a,1)}({\Bbb R}^n \times I), \text{ supported for }x\in\overline\Omega .
$$
 
When $p>1/a$,  $H_p^{a(2a)}(\overline{\Bbb
  R}^n_+ )$ consists, as recalled earlier in (2.21)ff., of the functions $
v=w+x_n^aK_0\varphi$ with $ w\in \dot H_p^{2a}(\overline{\Bbb  R}^n_+
)$, $ \varphi \in B_p^{a-1/p}({\Bbb  R}^{n-1})$; here $K_0$ is the
  Poisson operator $K_0\colon\varphi  \mapsto
e^{-\langle{D'}\rangle x_n}\varphi $ solving the Dirichlet problem for
  $1-\Delta $ with $(1-\Delta )u=0$ and nontrivial boundary data $\varphi $. Also in the curved case,
  $H_p^{a(2a)}(\overline\Omega )$ consists of $\dot
H_p^{2a}(\overline\Omega )$ plus a space of Poisson solutions in $\ol H_p^{a}(\Omega )$ multiplied
  by $d^a$.
\endexample

We observe moreover that if we take $r\ge 0$ such that
$$
r\;\cases =2a\text{ if }0<a<1/p,\\ <a+1/p \text{ if }1/p\le a<1 ,\endcases\tag4.36
$$
(note that $r\le 2a$ in all cases), then $
H_p^{a(2a)}(\comega)\subset \dot H_p^r(\comega)$,
for all $0<a<1$ (in view of (2.19), cf.\ also (4.21)). The statement (4.34) then implies that 
$$
u\in \ol H_p^{(r,r/(2a))}(\Omega \stimes I).\tag4.37
$$
When $f$ has a higher regularity than $L_p(\Omega \times I )$,  the
{\it interior} regularity can be improved a little  by use of Theorem 4.4 (for
the boundary regularity, see Theorem 4.18--19):

\proclaim{Theorem 4.12} Let $u$ be as in  Theorem {\rm 4.9}, and let
$r$ satisfy {\rm (4.36)}. Then $u$ satisfies {\rm (4.37)}, and
moreover, for $0<s\le r$,
$$
f\in  \ol H_p^{(s,s/(2a))}(\Omega \stimes I )\implies u\in H_{p,\operatorname{loc}}^{(s+2a,s/(2a)+1)}(\Omega \stimes I ).
\tag4.38
$$

In particular, if $a<1/p$ and $f\in  \ol H_p^{(2a,1)}(\Omega \stimes I
)$, then $u\in H_{p,\operatorname{loc}}^{(4a,2)}(\Omega \stimes I).$

\endproclaim

\subhead 4.3 Higher time-regularity \endsubhead
 
The result of Theorem 4.9 can be considerably extended by use of the
theory of Amann \cite{A97}, in the question of time-regularity. Fix $p\in\,]1,\infty [\,$. The fact that 
$-P_{\operatorname{Dir},p}$ is the generator of a
bounded holomorphic semigroup $T_p(t)$ in $L_p(\Omega )$ (for $t$ in a
sector around $\rp$ depending on $p$, cf.\ \cite{D89} Th.\ 1.4.2), assures that
there is an obtuse sector 
$$
V_\delta =\{\lambda \in{\Bbb C}\mid \arg\lambda \in \,]\pi /2-\delta ,3\pi /2+\delta [\,\},\tag4.39
$$
where the resolvent  $(P_{\operatorname{Dir},p}-\lambda )^{-1}$ exists
and  
satisfies an inequality
$$
|\lambda |\|(P_{\operatorname{Dir},p}-\lambda
)^{-1}\|_{\Cal L(L_p)}\le M,\tag4.40
$$
cf.\ Hille and Phillips \cite{HP57} Th.\ 17.5.1, or e.g.\ Kato \cite{K66}
Th. IX.1.23. Since   $P_{\operatorname{Dir},p}$ has a bounded
inverse, $\|(P_{\operatorname{Dir},p}-\lambda
)^{-1}\|_{\Cal L(L_p)}\le c$ also holds for $\lambda $ in a
neighborhood of 0, so we can replace $|\lambda |$ by $\ang\lambda $ in
(4.40) (with a larger constant $M'$). Note that furthermore, 
$$
\|P_{\operatorname{Dir},p}(P_{\operatorname{Dir},p}-\lambda
)^{-1}f\|_{L_p}=\|f+\lambda (P_{\operatorname{Dir},p}-\lambda
)^{-1}f\|_{L_p}\le (1+M)\|f\|_{L_p},
$$
so that $(P_{\operatorname{Dir},p}-\lambda )^{-1}$ is bounded
uniformly in $\lambda $ from $L_p(\Omega )$ to
$D(P_{\operatorname{Dir},p})$ with the graph-norm.

Thus, if we set
$$
E_0=L_p(\Omega ),\quad E_1=D(P_{\operatorname{Dir},p})=H_p^{a(2a)}(\comega),\tag4.41
$$
we have that $A=P_{\operatorname{Dir},p}$ satisfies (with
$B_\varepsilon =\{|\lambda |<\varepsilon \}$ for a small $\varepsilon >0$)
$$
\ang\lambda \|(A-\lambda
)^{-1}\|_{\Cal L(E_0)}+
\|(A-\lambda )^{-1}\|_{\Cal L(E_0,E_1)}\le c
\text{ for }\lambda \in V_\delta \cup B_\varepsilon .\tag4.42
$$
We can then apply Theorem 8.8 of Amann \cite{A97}. It is formulated with
vector-valued Besov spaces $B^s_{q,r}({\Bbb R};X)$ (valued in a Banach space
$X$, a function space in the applications), where the case $q=r=\infty $ is particularly interesting for our
purposes, since $B^s_{\infty ,\infty }({\Bbb R};X)$ equals the
vector-valued H\"older-Zygmund space $ C^s_*({\Bbb R};X)$. This
coincides with the vector-valued H\"older space $ C^s({\Bbb R};X)$
when $s\in \rp\setminus {\Bbb N}$, see  Remark 2.7 for further information.

\proclaim{Theorem 4.13}{\rm (From \cite{A97} Th.\ 8.8.)} Let
$s\in{\Bbb R}$ and $q,r\in [1,\infty ]$. Let $E_0$
and $E_1$ be Banach spaces, with $E_1$ continuously injected in $E_0$,
and let $A$ be a linear operator in $E_0$ with domain $E_1$ and range $E_0$,
satisfying {\rm (4.42)}.

For $f$ given in $B^s_{q,r}({\Bbb
R};E_0)\cap L_{1,\operatorname{loc}}({\Bbb R}; E_0)$ and supported in
$\crp$, the problem $$
Au+\partial_tu=f\text{ for }t\in {\Bbb R},\quad \supp u\subset \crp,
$$
 has a unique solution $u$ in the space of distributions in $  B^{s+1}_{q,r}({\Bbb
R};E_0)\cap B^{s}_{q,r}({\Bbb
R};E_1)$ supported in
$\crp$; it is described by
$$
u(t)=\int_0^t e^{-(t-\tau )A}f(\tau )\, d\tau , \quad t>0.
$$

\endproclaim

An application of this theorem 
with $q=r=\infty $ and $E_0,E_1$ defined
by (4.41) gives:

\proclaim{Theorem 4.14} Assumptions as in Theorem {\rm 4.9}. The
solution $u$ satisfies: When
$f\in \dot C^s_*(\crp ;L_p(\Omega ))\cap L_{1,\operatorname{loc}}({\Bbb R},L_p(\Omega ))$ for some $s\in{\Bbb R}$, then
$$
u\in \dot C^s_*(\crp ;H_p^{a(2a)}(\comega))\cap \dot C^{s+1}_*(\crp ;L_p(\Omega )).\tag4.43
$$

In particular, for $s\in \rp\setminus{\Bbb N}$,
$$
f\in \dot C^s(\crp ;L_p(\Omega ))\iff u\in \dot
C^s(\crp ;H_p^{a(2a)}(\comega))\cap \dot C^{s+1}(\crp ;L_p(\Omega )). \tag 4.44
$$ 
\endproclaim

We here use conventions as in (2.4): The set of 
$u(x,t)\in C^s_*({\Bbb R};X)$ that vanish for $t<0$ is
denoted  $\dot C^s_*(\crp ;X)$. The conclusion $\implies$ is a special
case of Theorem 4.13 with $A=P_{\operatorname{Dir},p}$; the converse
$\impliedby$ follows immediately by application of $A+\partial _t$ to $u$.

\example{Remark 4.15} Amann's theorem can of course also be
applied with other choices of $B^s_{q,r}$, e.g.\ 
Sobolev-\-Slobodetski\u\i{} spaces $W^{s,q}=B^s_{q,q}$ when $s\in \rp\setminus {\Bbb N}$;
this leads to analogous statements. (Note that integer cases $s\in{\Bbb
N}_0$ are not covered, in particular not the result of Theorem 4.9 for
$p\ne 2$, since 
the $H^s_p$-spaces are not in the scale $B^s_{q,r}$ when $p\ne 2$.)
\endexample

Observe the consequences:

\proclaim{Corollary 4.16} Assumptions as in Theorem {\rm 4.9}. The
solution $u$ satisfies: 
$$f\in \dot C^\infty (\crp ;L_p(\Omega ))\iff
u\in \dot C^\infty (\crp ;H_p^{a(2a)}(\comega)).\tag4.45
$$
Moreover,
$$
\aligned
f\in \dot C^s(\crp ;L_\infty (\Omega ))&\implies
u\in \dot C^{s} (\crp ;d^a\ol C^{a-\varepsilon }(\Omega)),\\
f\in \dot C^\infty (\crp ;L_\infty (\Omega ))&\implies
u\in \dot C^\infty (\crp ;d^a\ol C^{a-\varepsilon }(\Omega)),
\endaligned
\tag4.46
$$
for $s\in \rp\setminus {\Bbb N}$ and   $\varepsilon \in \,]0,a]$.
\endproclaim

\demo{Proof} When $f\in \dot C^\infty (\crp ;L_p(\Omega ))$, it is in
particular in $ \dot C^s (\crp ;L_p(\Omega ))$ for all
$s\in\rp\setminus {\Bbb N}$, where (4.44) implies that $ u\in \dot
C^s(\crp ;H_p^{a(2a)}(\comega))$. Taking intersections over $s$, we
find (4.45).

(4.46) follows from (4.44) resp.\ (4.45), since $L_\infty (\Omega )\subset L_p(\Omega )$ for
$1<p<\infty $, and $\bigcap_pH_p^{a(2a)}(\comega)\subset d^a\ol
C^{a-\varepsilon }(\Omega )$ (as already observed in \cite{G15}
Sect.\ 7). \qed
\enddemo  

The statements in (4.46) improve the results of \cite{FR17} and
\cite{RV18} concerning time-regularity.

Amann's theorem  can also be used in situations with a  higher
$x$-regularity of $f$. The crucial point is to obtain an estimate
(4.42) on the desired pair of Banach-spaces $E_0,E_1$ with $E_1\subset E_0$.
 
\proclaim{Lemma 4.17} Let $1<p<\infty $, and let $A$ be an operator
acting like $P_{\operatorname{Dir},p}$ in Theorem {\rm 4.9}. 
 When
$2a<1/p$, $A$ is bijective from $H^{a(4a)}(\comega)$ to $\dot
H^{2a}_p(\comega)$, and we have a
resolvent inequality for $\lambda \in V_\delta \cup B_\varepsilon $ as
in {\rm
 (4.42)}:
$$
\ang\lambda \|(A-\lambda )^{-1}f\|_{\dot H_p^{2a}}+\|(A-\lambda
)^{-1}f\|_{ H_p^{a(4a)}}\le c \|f\|_{\dot H^{2a}_p}, \text{ for }f\in
\dot H_p^{2a}(\comega).\tag 4.47
$$

\endproclaim

\demo{Proof}
Define $D_{s,p}(A)$ by 
$$
D_{s,p}(A)=H_p^{a(2a+s)}(\comega);\tag4.48
$$
it will be used for $s=0$ and $2a$. 
Since $2a<1/p$, $\ol H^{2a}_p(\Omega )$ identifies
with $\dot H^{2a}_p(\comega)$, so $A$ has the  asserted bijectiveness property
in view of Theorem 2.5 plus  \cite{G14} Th.\ 3.5
on the invariance of kernel and cokernel. 

By (4.42), we have the inequality 
$$
\ang\lambda \|(A-\lambda )^{-1}f\|_{L_p}+\|(A-\lambda
)^{-1}f\|_{D_{0,p}}\le c\|f\|_{L_p}\text{ for }f\in L_p(\Omega ),\tag4.49
$$
and we want to lift it to the case where $L_p(\Omega )$ is replaced by
$E_0=\dot
H_p^{2a}(\comega)$ and $D_{0,p}(A)$ is replaced by  $E_1=D_{2a,p}(A)$.
Let $f\in \dot H_p^{2a}(\comega)=D_{0,p}(A)$, and denote $Af=g$; it
lies in $L_p(\Omega )$, and $\|g\|_{L_p}\simeq \|f\|_{\dot H_p^{2a}}$. Then 
$$
\aligned
\|\lambda (A-\lambda )^{-1}f\|_{\dot H_p^{2a}}&=\|-f+A(A-\lambda
)^{-1}f\|_{\dot H_p^{2a}}=\|-f+(A-\lambda )^{-1}Af\|_{\dot
H_p^{2a}}\\
&\le \|f\|_{\dot H_p^{2a}}+\|(A-\lambda )^{-1}g\|_{\dot
H_p^{2a}}\le \|f\|_{\dot H_p^{2a}}+c\|g\|_{L_p}\le c'\|f\|_{\dot H_p^{2a}},
\endaligned\tag4.50
$$
where we applied (4.49) to $g$.
Next, 
$$
\aligned
\|(A-\lambda )^{-1}f\|_{D_{2a,p}(A)}&\simeq \|A(A-\lambda
)^{-1}f\|_{\ol H_p^{2a}}\simeq \|A(A-\lambda
)^{-1}f\|_{\dot H_p^{2a}}\\
&= \|f+\lambda (A-\lambda
)^{-1}f\|_{\dot H_p^{2a}}\le (1+c')\|f\|_{\dot H_p^{2a}},
\endaligned
$$
using (4.50). Together with (4.50), this proves (4.47).\qed
\enddemo

This leads to a supplement to Theorem 4.14:

\proclaim{Theorem 4.18} Assumptions as in Theorem {\rm 4.9}. 
If $2a<1/p$, the
solution satisfies,
for $s\in \rp\setminus{\Bbb N}$:
$$
\aligned
f\in \dot C^s(\crp ;\dot H_p^{2a}(\comega ))&\iff u\in \dot
C^s(\crp ;H_p^{a(4a)}(\comega))\cap \dot C^{s+1}(\crp ;\dot H_p^{2a}(\comega )),\\
f\in \dot C^\infty (\crp ;\dot H_p^{2a}(\comega ))&\iff u\in \dot
C^\infty (\crp ;H_p^{a(4a)}(\comega)). 
\endaligned \tag 4.51
$$ 
\endproclaim

\demo{Proof} We apply Theorem 4.13 with $E_0=\dot
H^{2a}(\comega)$, $E_1=H^{a(4a)}(\comega)$, where the required
resolvent inequality is shown in Lemma 4.17. \qed
\enddemo

Let us also apply Amann's theorem to the $L_2$-operators studied in
Theorem 4.6, listing
just the resulting H\"older estimates. Here when $I=\,]0,T[\,$, we
denote by $C^s_+(I;X)$ the space of functions in $\ol C^s(\,]-\infty
,T[\,;X)$ with support in $[0,T]$.

\proclaim{Theorem 4.19} Assumptions as in Theorem {\rm 4.6}. 
With $\xi _0$ as in {\rm(4.18)}, there is a set $V_\delta \cup B_\varepsilon $
(cf.\ {\rm (4.39)}) such that $A=P_{\operatorname{Dir},2}+\xi _0$
satisfies an
inequality {\rm (4.42)}, both for the choice $\{E_0, E_1\}=\{L_2(\Omega ),
H^{a(2a)}(\comega)\}$ when $a<1$, and for the choice $\{E_0,E_1\}=\{\dot
H^{2a}(\comega),H^{a(4a)}(\comega)\}$ when $a<\frac14$.

The
solution of {\rm (4.15)} satisfies,
for $s\in \rp\setminus{\Bbb N}$, $I=\,]0,T[\,$:
$$
\align
f\in  C^s_+(I ;L_2(\Omega ))&\iff u\in 
C^s_+(I ;H^{a(2a)}(\comega))\cap  C^{s+1}_+(I ;L_2(\Omega
)),\tag4.52\\
f\in  C^\infty _+(I ;L_2(\Omega ))&\iff u\in 
C^\infty _+(I ;H^{a(2a)}(\comega)).\tag4.53
\endalign
$$
Moreover, if $a<\frac14$,
$$
\aligned
f\in C^s_+(I ;\dot H^{2a}(\comega ))&\iff u\in 
C^s_+(I ;H^{a(4a)}(\comega))\cap  C^{s+1}_+(I ;\dot H^{2a}(\comega )),\\
f\in  C^\infty _+(I ;\dot H^{2a}(\comega ))&\iff u\in 
C^\infty _+(I ;H^{a(4a)}(\comega)). 
\endaligned \tag 4.54
$$ 
\endproclaim

\demo{Proof} 
The estimate of $(A-\lambda )^{-1}$ with $E_0=L_2(\Omega )$,
$E_1=H^{a(2a)}(\comega)$, is assured by the information in the proof
of Theorem 4.6 that the numerical range and spectrum of
$P_{\operatorname{Dir},2}$ is contained in an angular set $\{z\in{\Bbb C}\mid
|\operatorname{Im}z|\le C(\operatorname{Re}z+\xi ), \,
\operatorname{Re}z>\xi _0\}$. (This is a standard fact in the
theory of operators defined from sesquilinear forms, see e.g.\
Cor.\ 12.21 in \cite{G09}.) Then the resolvent estimate holds for
$\lambda \to\infty $ on the rays in a closed sector disjoint from the sector $\{
|\operatorname{Im}z|\le C\operatorname{Re}z \}$. 

The estimate of $(A-\lambda )^{-1}$ with $E_0=\dot H^{2a}(\comega )$,
$E_1=H^{a(4a)}(\comega)$ now follows exactly as in Lemma 4.17, when $a<\frac14$.

For the solvability assertions, note that $u(x,t)$ satisfies
$P_{\operatorname{Dir},2}u+\partial_tu=f$ if and only if
$v(x,t)=e^{-\xi _0t}u(x,t)$ satisfies
$(P_{\operatorname{Dir},2}+\xi _0)v+\partial_tv=e^{-\xi _0t}f$. 
 
An application of Theorem 4.13 to $A=P_{\operatorname{Dir},2}+\xi _0$ leads to a solvability result
like (4.44) for a solution $v$ with right-hand
side $f_1$. When the problem (4.15) is considered for a given function $f\in
C^s_+(I; E_0)$, we extend $f$ to a function $\tilde f\in \dot
C^s(\crp;E_0)$, and apply the solvability result with $f_1=e^{-\xi _0t}\tilde f$;
this gives a solution $v$, and we set $u=(e^{\xi _0t}v)|_{I}$, solving
(4.15). It has the regularity claimed  in (4.52), and conversely,
application of 
$P_{\operatorname{Dir},2}+\partial_t$ to such a function gives a
right-hand side in $C^s_+(I; E_0)$. (4.53) follows immediately. 

The statements in (4.54) follow by a similar application of Theorem
4.13 with $E_0=\dot H^{2a}(\comega )$,
$E_1=H^{a(4a)}(\comega)$. \qed 

\enddemo

It seems plausible that (4.51) and (4.54) can be extended 
to all $a<1$ when
$\dot H_p^{2a}(\comega)$ is replaced by the possibly weaker space
$\dot H_p^{r}(\comega)$, cf.\ (4.36), and $H_p^{a(4a)}(\comega)$ is
replaced by $H_p^{a(2a+r)}(\comega)$ (for $p=2$, it is to some extent
obtained in (4.26)). This might be based on an extension of Amann's
strategy. The
Fourier transform used in \cite{LM68} is in \cite{A97}
replaced by multiplier theorems in a 
suitable sense; it could be investigated whether something similar
might be done based on resolvent inequalities in the
style of (4.32), (4.33).

On the other hand, we expect that the limitations on $r$ are essential in
some sense, since the domain and range for the Dirichlet problem are
not compatible in higher-order spaces, as noted also earlier.

\example{Remark 4.20}
The results in this subsection on higher $t$-regularity based on Theorem 4.13 from \cite{A97}, 
are in the case $p\ne 2$
restricted to the $x$-independent case $P=\Op(p(\xi ))$ with
homogeneous symbol. However, since the proofs rely entirely on general
resolvent inequalities, there is a leeway to extend the results to
perturbations $P+P'$ when $P'$ is so small that such
resolvent inequalities still hold; this will allow some $x$-dependence and
lower-order terms. The idea can be further developed.
\endexample

It is very likely that other methods could be useful and bring out
further perspectives  for the heat problem. One possibility could be to
try to establish an $H^\infty $-calculus for the 
realizations of fractional Laplacians and their generalizations.
This was done for  elliptic differential
operators with boundary conditions in a series of works, see e.g.\ Denk, Hieber and Pr\"uss 
\cite{DHP03} for
an explanation of the theory and many references. However,
pseudodifferential problems pose additional difficulties. 
The only contribution treating an operator within the Boutet de
Monvel calculus, that we know of, is Abels
\cite{A05a}. The present pseudodifferential operators do not even belong to the Boutet de
Monvel calculus, but have a different boundary behavior.

It will be interesting to see to what extent the results can be further
improved by this or other tools from the extensive literature on differential heat
operator problems.

\end{document}

%% file: gmacro2.tex
\def\supp{\operatorname{supp}}

\def\crp{\overline{\Bbb R}_+}

\def\rnp{{\Bbb R}^n_+}
\def\rnm{\Bbb R^n_-}

\def\crnp{\overline{\Bbb R}^n_+}

\def\comega{\overline\Omega }

\def\Rn{\Bbb R^n}
\def\ang#1{\langle {#1} \rangle}

\def\Op{\operatorname{Op}}
\def\Pfrac{\tsize\frac1{\raise 1pt\hbox{$\scriptstyle p$}}}
\def\pfrac{\frac1{\raise 1pt\hbox{$\scriptscriptstyle p$}}}
\def\Pfracc#1{\tsize\frac{#1}{\raise 1pt\hbox{$\scriptstyle p$}}}
\def\pfracc#1{\frac{#1}{\raise 1pt\hbox{$\scriptscriptstyle p$}}}

\def\Zfrac{\tsize\frac1{\raise 1pt\hbox{$\scriptstyle z$}}}
\def\zfrac{\frac1{\raise 1pt\hbox{$\scriptscriptstyle z$}}}

\def\rp{ \Bbb R_+}

\define\stimes{\!\times\!}

\def\R{\Bbb R}

\def\ol{\overline}
\def\SD{\Cal S}
\def\E{\Cal E}
\def\F{\Cal F}

%% file: isaac9.bbl
\Refs
\widestnumber\key{[CWWW18]}

\ref\no[A15] \by N.\ Abatangelo \paper Large s-harmonic functions and
boundary blow-up solutions for the fractional Laplacian \jour
 Discrete Contin.\ Dyn.\ Syst.\ \vol 35 \yr2015\pages 5555--5607
\endref 

\ref\no[AJS18] \by N. Abatangelo, S. Jarohs and A. Saldana \paper 
Integral representation of solutions to higher-order
fractional Dirichlet problems on balls \jour
  Comm. Contemp. Math. \finalinfo to appear, arXiv: 1707.03603 
\endref 


\ref\no[A05] \by H.\ Abels \paper Pseudodifferential boundary value
problems with non-smooth coefficients
\jour Comm. Part. Diff. Equ.\vol 30 \yr 2005 \pages 1463-–1503\endref 

\ref\no[A05a] \by H.\ Abels \paper Reduced and generalized Stokes
resolvent equations in asymptotically flat layers. II. $H_\infty $-calculus
\jour J. Math. Fluid Mech.\vol 7 \yr 2005 \pages 223-–260\endref 

\ref\no[A97] \by H. Amann \paper  Operator-valued Fourier multipliers, vector-valued Besov spaces, and
 applications\jour Math. Nachr. \vol 186 \yr 1997 \pages 5--56\endref


\ref\no[A04]\by D. Applebaum \paper L\'evy processes --— from probability to
finance and quantum groups\jour Notices Amer. Math. Soc. \vol 51
\yr2004 \pages 1336-–1347 \endref

\ref\no[BWZ17]\by U. Biccari, M. Warma and E. Zuazua \paper Local
elliptic regularity for the Dirichlet fractional Laplacian \jour
Advanced Nonlinear Studies \vol 17 \yr 2017 \pages 387--409\endref

\ref\no[BWZ17a]\by U. Biccari, M. Warma and E. Zuazua \paper Addendum: Local
elliptic regularity for the Dirichlet fractional Laplacian \jour
Advanced Nonlinear Studies \vol 17 \yr 2017 \pages 837
\endref


\ref\no[BWZ18]\by U. Biccari, M. Warma and E. Zuazua \paper Local
regularity for fractional heat equations \finalinfo arXiv: 1704.07562
\endref

\ref\no[BBC03]\paper     Censored stable processes
 \by   K. Bogdan, K. Burdzy and Z.-Q. Chen \vol 127\pages 89--152
\jour Prob. Theory Related Fields\yr 2003 
\endref


\ref\no[BHL16]\by T. Boulenger, D. Himmelsbach and E. Lenzmann \paper
Blowup for fractional NLS\jour J. Funct. Anal.\vol 271 \yr 2016 
\pages 2569-–2603
\endref 

\ref\no[B71] \by L.\ Boutet de Monvel\paper Boundary problems for pseudo-differential
operators \jour Acta Math.\ \vol 126 \yr1971 \pages 11--51\endref
 
\ref\no[CS07] \by L.\ Caffarelli and  L.\  Silvestre\paper An extension problem related
to the fractional Laplacian \jour Comm.\ Part.\ Diff.\ Eq.\ \vol 32 \yr
2007 \pages 1245--1260
\endref

\ref\no[CSS08] \by   L. A. Caffarelli, S. Salsa and Luis Silvestre
\jour Invent. Math.
\yr 2008 \vol 171 \pages 425-–461
\paper Regularity estimates for the solution and the free boundary of the obstacle problem for the fractional Laplacian 
\endref

\ref\no[CD14]
\by H. Chang-Lara and G. Davila \paper Regularity for solutions of non
local parabolic equations\jour Calc. Var. Part. Diff.
Equations \vol 49 \yr2014 \pages 139-–172
\endref

\ref\no[CS98]\by Z.-Q. Chen and R. Song \paper Estimates on Green
functions and Poisson kernels for symmetric stable processes \jour
Math. Ann. \vol 312 \yr1998 \pages 465--501
\endref 

\ref\no[CWWW18] \by M. Chen,  B. Wang,  S. Wang and  M. W. Wong 
\paper
On dissipative nonlinear evolutional pseudo-differential
equations
\finalinfo arXiv:1708.09519
\endref


\ref\no[CT04] \by R. Cont and P. Tankov \book Financial modelling with
jump processes, Chapman \& Hall/CRC Financial Mathematics Series \publ
Chapman \& Hall/CRC \publaddr Boca Raton, FL \yr 2004
\endref

\ref\no[D89] \by E. B. Davies\book Heat kernels and spectral
theory. Cambridge Tracts in Mathematics, 92\publ Cambridge University
Press \publaddr Cambridge \yr 1989 
\endref

\ref\no[DHP03] \by R. Denk, M. Hieber and J. Pr\"uss \paper R-boundedness, Fourier
 multipliers and problems of elliptic and parabolic
 type \jour Mem. Amer. Math. Soc. \vol 166 \yr 2003\finalinfo no. 788, viii+114
 pp \endref 

 \ref\no[DG17] \by S. Dipierro and H. Grunau \paper Boggio's formula
 for fractional polyharmonic Dirichlet problems \jour Ann. Mat. Pura
 Appl.\vol 196 \yr2017 \pages 1327-–1344 \endref

\ref\no[DKK17] \by
B.  Dyda, A. Kuznetsov and M. Kwasnicki \paper Eigenvalues
 of the fractional Laplace operator in the unit
 ball \jour J. Lond. Math. Soc. (2)\vol 95 \yr2017 \pages 500-–518 \endref 

\ref\no[E81] \by G.\ Eskin\book Boundary value problems for elliptic
pseudodifferential equations,  AMS Translations \yr1981 \publ
Amer. Math. Soc. \publaddr Providence, R. I.\endref


\ref\no[FK13]
\by M. Felsinger and M. Kassmann \paper Local regularity for parabolic
nonlocal operators\jour Comm. Part. Diff. Equations \vol 38
\yr2013 
\pages 1539-–1573\endref

\ref\no[FR17] 
\paper Regularity theory for general stable operators: parabolic equations
\by
X. Fernandez-Real and X. Ros-Oton\jour
J. Funct. Anal. \vol 272 \yr2017 \pages 4165--4221
\endref


\ref\no[FG16]\by R. Frank and L. Geisinger \paper Refined
semiclassical asymptotics for fractional powers of the Laplace
operator\jour J. Reine Angew. Math. \vol712 \yr 2016 \pages 1-–37
\endref


\ref\no[FOT94]\by M. Fukushima, Y. Oshima and M. Takeda \book
Dirichlet forms and symmetric Markov processes. De Gruyter Studies in
Mathematics, 19\publ Walter de Gruyter \& Co.\publaddr Berlin \yr 1994
\endref


\ref \no[GMS12]\by M. Gonzalez, R. Mazzeo and Y. Sire \paper Singular
solutions of fractional order conformal Laplacians \jour
J. Geom. Anal.\vol 22 \yr2012\pages 845-–863 \endref

\ref\no[G90] \by G.\ Grubb\paper Pseudo-differential boundary problems in Lp spaces,
Comm.\ Part.\ Diff.\ Eq.\ \vol 15 \yr1990 \pages 289--340\endref 

\ref\no[G95]
\by G. Grubb\paper Parameter-elliptic and parabolic pseudodifferential
boundary problems in global Lp Sobolev spaces \jour Math. Z. \vol 218
\yr 1995 \pages 43--90
\endref

 \ref\no[G96]\by 
{G.~Grubb}\book Functional calculus of pseudodifferential
     boundary problems.
 Pro\-gress in Math.\ vol.\ 65, Second Edition \publ  Birkh\"auser
\publaddr  Boston \yr 1996
\endref

\ref\no[G09]\by G. Grubb\book Distributions and operators. Graduate
Texts in Mathematics, 252 \publ Springer \publaddr New York\yr 2009
 \endref


\ref\no[G14] \by G.\ Grubb\paper  
Local and nonlocal boundary conditions for $\mu $-transmission
and fractional elliptic pseudodifferential operators \jour 
 Analysis and P.D.E.\  \vol 7 \yr 2014\pages 1649--1682\endref

\ref\no[G15] \by G.\ Grubb\paper Fractional Laplacians on domains, 
a development of H\"o{}rmander's theory of $\mu$-transmission
pseudodifferential operators \jour Adv.\ Math.\  \vol 268 \yr 2015 \pages
478--528\endref

\ref\no[G15a] \by G.\ Grubb\paper Spectral results for mixed problems
and fractional elliptic operators \jour J. Math. Anal. Appl.  \vol 421 \yr 2015 \pages
1616--1634\endref

\ref\no[G16] \by  G.\ Grubb \paper Regularity of spectral
fractional Dirichlet and Neumann problems \jour Math.\
Nachr.\  \vol 289 \yr 2016\pages 831--844
\endref


\ref\no[G16a] \by  G.\ Grubb \paper Integration by parts and  Pohozaev
identities for space-dependent fractional-order operators \jour J.\
Diff.\ Eq.\ \vol 261 \yr 2016\pages 1835--1879
\endref

\ref\no[G18]\by G. Grubb \paper  Regularity in $L_p$ Sobolev spaces
of solutions to fractional heat equations \jour J. Funct. Anal.\finalinfo 
published online  https://doi.org/10.1016/j.jfa.2017.12.011  (arXiv:1706.06058)
\endref

\ref\no[G18a]\by G. Grubb \paper   Green's formula and a Dirichlet-to-Neumann operator for
fractional-order pseudodifferential operators \finalinfo 
arXiv:1611.03024, to appear in Comm. Part. diff. Eq.
\endref

\ref\no[G18b]\by G. Grubb \paper   Limited regularity of solutions to
fractional heat equations \finalinfo (in preparation)
\endref

\ref
\no  [GS90] 
\by G. Grubb and V. A. Solonnikov
\paper Solution of parabolic pseu\-do-\-dif\-fe\-ren\-ti\-al in\-i\-ti\-al-boun\-d\-a\-ry
value problems
\jour J. Diff. Equ. 
\vol 87
\yr1990
\pages 256--304
\endref

\ref \no[HP57]\by E. Hille and R. S. Phillips \book Functional
analysis and semi-groups. rev. ed. American Mathematical Society
Colloquium Publications, vol. 31 \publ American Mathematical Society
\publaddr Providence, R. I. 
\yr 1957\endref

\ref\no[H65] \by L.\ H\"o{}rmander\paper Ch.\ II, Boundary problems for
``classical'' pseudo-differential operators. Unpublished
lecture notes at Inst.\ Adv.\ Studies, Princeton 1965 \finalinfo \TeX-typed version available at
 http://www.math.ku.dk$\sim$grubb/LH65.pdf
\endref


\ref\no[H85] \by L.\ H\"ormander\book The analysis of linear partial
differential operators, III \publ Springer Verlag \yr 1985 \publaddr
Berlin
\endref

\ref\no[J02] \by T. Jakubowski \paper The estimates for the Green function in Lipschitz
 domains for the symmetric stable processes \jour
 Probab. Math. Statist. \vol 22
 \yr 2002 \pages  419-441\endref

\ref\no[JX15]\by 
T. Jin and J. Xiong \paper Schauder estimates for solutions of linear
parabolic integro-differential equations \jour Discrete
Contin. Dyn. Syst. \vol 35 \yr2015 \pages 5977-–5998\endref

\ref \no [K66]\by T. Kato  \book  Perturbation theory for linear
operators. Die Grundlehren der mathematischen Wissenschaften, Band
132\publ Springer-Verlag New York, Inc. \publaddr New York 
\yr 1966\endref 

\ref\no[K97] \by T. Kulczycki \paper Properties of Green function of symmetric stable
 processes \jour Probab. Math. Statist. \vol 17 \yr1997 \pages
 339-364 \endref

\ref\no[L87]\by D. Lamberton\paper \'E{}quations d'\'evolution
lin\'eaires associ\'ees \`a des semi-groupes de contractions dans les
espaces Lp \jour J. Funct. Anal. \vol 72 \yr 1987 \pages 252--262
\endref

\ref\no[LPPS15]\by T. Leonori, I. Peral, A. Primo and F. Soria \paper
Basic estimates for solutions of a class of nonlocal elliptic and
parabolic equations
\jour
Discrete Contin. Dyn. Syst. \vol 35 \yr 2015 \pages 6031--6068
\endref


\ref \no   [LM68] 
\by J.-L. Lions and E. Magenes 
\book Probl\`emes aux limites non homog\`enes et applications. Vol. 1
et 2
\yr1968 
\publ Editions Dunod 
\publaddr Paris
\endref

\ref \no[MNP18] \by F. Monard, R. Nickl and G. P. Paternain 
\paper
    Efficient nonparametric Bayesian inference for X-Ray transforms
     \finalinfo arXiv:1708.06332 
\endref 

\ref\no[R16]\paper
Nonlocal elliptic equations in bounded domains: a survey
\by X. Ros-Oton
\jour Publ. Mat. \vol60 \yr2016 \pages 3--26\endref

\ref\no[R18]\paper Boundary  regularity,  Pohozaev identities and
nonexistence results
\by X. Ros-Oton
\finalinfo \linebreak  arXiv:1705.05525, to appear as a chapter in "Recent Developments in the Nonlocal Theory" by De Gruyter\endref

\ref\no[RS14] \by X.\ Ros-Oton and J.\ Serra\paper The Dirichlet problem for the
fractional Laplacian: regularity up to the boundary \jour J.\ Math.\ Pures
Appl.\ \vol 101 \yr 2014 \pages  275--302\endref

\ref\no[RS14a] \by  X.\ Ros-Oton and J.\ Serra \paper
The Pohozaev identity for the fractional Laplacian \jour
Arch.\ Rat.\ Mech.\ Anal.\ \vol 213 \yr 2014 \pages 587--628\endref

\ref\no[RS15] \by  X.\ Ros-Oton and J.\ Serra \paper
Local integration by parts and Pohozaev identities for higher order
fractional Laplacians \jour  Discrete Contin. Dyn. Syst. \vol 35 \yr
2015 \pages 2131-–2150\endref

\ref\no[RSV17] \by X.\ Ros-Oton, J.\ Serra and E.\ Valdinoci\paper 
Pohozaev identities for anisotropic integro-dif\-fe\-ren\-tial operators
\jour Comm. Part. Diff. Eq. \vol 42 \yr 2017 \pages 1290--1321
\endref

\ref\no[RV18] \by X.\ Ros-Oton and H.\ Vivas \paper Higher-order
boundary regularity estimates for nonlocal parabolic equations
\finalinfo  arXiv:1711.02075 \endref

\ref
\no  [ST87] 
\by H.-J. Schmeisser and H. Triebel
\book Fourier analysis and function spaces
\publ Wiley
\publaddr New York
\yr1987
\endref

\ref \no[S01] \by E. Schrohe
\paper  A  short  introduction  to  Boutet  de  Monvel's  calculus
\inbook
Approaches
to Singular Analysis, Oper. Theory Adv. Appl. 125 \publ Birkh\"auser \publaddr
Basel \yr 2001
\pages 85-–116
\endref 

\ref\no[S69]\by R. Seeley \paper The resolvent of an elliptic boundary
problem\jour Amer. J. Math.\vol 91 \yr 1969 \pages 889-–920
\endref

\ref\no[T78] 
\by H. Triebel
\book Interpolation theory, function spaces, differential operators
\publ North-Hol\-land Publ. Co.
\publaddr Amsterdam, New York
\yr1978
\endref

\ref\no[ZG16] \by M. Zaba and P. Garbaczewski \paper Ultrarelativistic
bound states in the spherical well \jour J. Math. Phys. \vol 57
\yr2016 \pages 26 pp \endref

\endRefs
